\documentclass[a4,11pt]{article}
\usepackage{amsmath,amssymb,amstext}
\usepackage{graphicx}

\usepackage[tmargin=1.25in,bmargin=1.25in,rmargin=1.25in,lmargin=1.25in]{geometry}

\title{Successive Standardization of Rectangular Arrays}
\author{Richard A. Olshen\\ {\small Stanford University } \and Bala Rajaratnam\\ {\small Stanford University}}

\date{}

\newcommand{\Xb}{\mathbf X}

\begin{document}

\maketitle

\begin{abstract}
In this note we illustrate and  develop further with mathematics and examples, the work on successive standardization (or normalization) that is studied earlier by the same authors in \cite{olshen_rajaratnam1} and \cite{olshen_rajaratnam2}. Thus, we deal with successive iterations applied to rectangular arrays of numbers, where to avoid technical difficulties an array has at least three rows and at least three columns.  Without loss, an iteration begins with operations on columns:  first subtract the mean of each column; then divide by its standard deviation.  The iteration continues with the same two operations done successively for rows.  These four operations applied in sequence completes one iteration.  One then iterates again, and again, and again,.... In \cite{olshen_rajaratnam1} it was argued that if arrays are made up of real numbers, then the set for which convergence of these successive iterations fails has Lebesgue measure 0.  The limiting array has row and column means 0, row and column standard deviations 1. A basic result on convergence given in \cite{olshen_rajaratnam1} is true, though the argument in \cite{olshen_rajaratnam1} is faulty. The result is stated in the form of a theorem here, and the argument for the theorem is correct. Moreover, many graphics given in \cite{olshen_rajaratnam1} suggest that but for a set of entries of any array with  Lebesgue measure 0, convergence is very rapid, eventually exponentially fast in the number of iterations.  Because we learned this set of rules from Bradley Efron, we call it ``Efron's algorithm". More importantly, the rapidity of convergence is illustrated by numerical examples.
\end{abstract}

\noindent {\it Key words:} rectangular arrays, successive iterations, standardization, exponentially fast convergence.

\noindent {\it AMS 2000 subject classifications:} 62H05, 60F15, 60G46.

\section{Introduction}

Research summarized in this paper is an extension of that reported in \cite{olshen_rajaratnam1} and a conference proceeding \cite{olshen_rajaratnam2}, and concerns successive standardization (or normalization) of large rectangular arrays $\bf{X}$ of real numbers, such as arise in gene expression, from protein chips, or in the earth and environmental sciences. The basic message here is that convergence that holds on all but a set of measure 0 in the paper by Olshen and Rajaratnam (\cite{olshen_rajaratnam1} is shown here to be exponentially fast in a sense we make precise. The basic result in \cite{olshen_rajaratnam1}, though true is not argued correctly in \cite{olshen_rajaratnam1}. The gap is filled here. Typically there is one column per subject; rows correspond to ``genes" or perhaps gene fragments (including those that owe to different splicing of ``the same" genes, or proteins).  Typically, though not always, columns divide naturally into two groups: ``affected" or not. Two-sample testing of rows that correspond to ``affected" versus other individuals or samples is then carried out simultaneously for each row of the array.  Corrections for multiple comparisons may be very simple, or might perhaps allow for ``false discovery."  

As was noted in \cite{olshen_rajaratnam1}, data may still suffer from problems that have nothing to do with differences between groups of subjects or differences between ``genes" or groups of them.  There may be variation in background, perhaps also in ``primers." Thus, variability across subjects might be unrelated to status. In comparing two random vectors that may have been measured in different scales, one puts observations ``on the same footing" by subtracting each vector's mean and dividing by its standard deviation.  Thereby, empirical covariances are changed to empirical correlations, and comparisons proceed.  But how does one do this in the rectangular arrays described earlier?  An algorithm by which such ``regularization" is accomplished was described to us by colleague Bradley Efron; so we call the full algorithm \emph{Efron's algorithm}. We shall use the terms successive ``normalization" or successive ``standardization" interchangeably (and also point out that Bradley Efron considers the the latter term a better description of the algorithm). To avoid technical problems that are described in \cite{olshen_rajaratnam1} and repeated here, we assume there are at least three rows and at least three columns to the array.  It is immaterial to convergence, though not to limiting values, whether we begin regularization to be described by row or by column.  In order that we fix an algorithm, we begin by column in the computations though by row in the mathematics. Thus, we first mean polish the column, then standard deviation polish the column; next we mean polish the row, and standard deviation polish the row.  The process is then repeated.  By ``mean polish" of, say, a column, we mean subtract the mean value for that column from every entry.  By ``standard deviation" polish the column, we mean divide each number by the standard deviation of the numbers in that column.  Definitions for ``row" are entirely analogous.

In \cite{olshen_rajaratnam1} convergence is studied with entries in the rectangular array with I rows and J columns viewed as elements of  $\mathbb{R}^{IJ}$, Euclidean IJ space. Convergence holds for all but a Borel subset of  $\mathbb{R}^{IJ}$ of Lebesgue measure 0.  The limiting vector has all row and column means 0, all row and column standard deviations 1. We emphasize that to show convergence on a Lebesgue set of full measure, it is enough to find a probability mutually absolutely continuous with respect to Lebesgue measure for which convergence is established with probability 1.

%

\section{Preliminaries}

We begin precise definitions:

$$ {\bar X}^{(0)}_{i \cdot}  = {1\over k}\ \sum^k_{i=1}\ X_{ij}\ $$

$$(S^{(0)}_i)^2 = {1\over k}\ \sum^k_{j=1}\ (X_{ij} - {\bar X}^{(0)}_{i \cdot})^2 = ({1\over k}\ \sum^k_{j=1} \ X^2_{ij}  - k ({\bar X}^{(0)}_{i \cdot})^2  ) $$
 
$${\bf X}^{(1)} = \{ X_{ij}^{(1)} \},\; \textrm{where} \; X^{(1)}_{ij} = {X_{ij} - {\bar X}^{(0)}_{i \cdot} \over
{S^{(0)}_i}}$$
 
By analogy, set
 
$${\bf X}^{(2)} = \{ X_{ij}^{(2)} \},$$ 
where
$$X^{(2)}_{ij} = {X^{(1)}_{ij} - {\bar X}^{(1)}_{\cdot j} \over {S^{(1)}_j} } \; \textrm{and} \;{\bar X}^{(1)}_{\cdot j} = {1\over n} \ \sum^n_{i=1} \ X^{(1)}_{ij}$$ 

Now,
 
$$(S^{(1)}_j)^2 = {1\over n} \ \sum^n_{i=1} (X^{(1)}_{ij} - {\bar X}^{(1)}_{\cdot j})^2,$$

and analogously define $(S^{(n)}_j)^2$ for $n = 2,3,..$ By induction on dimension applied to regular conditional distributions, $(S^{(0)}_i)^2 > 0$ a.s.
 
As explained in \cite{olshen_rajaratnam1} the $2 \times 2$ case illustrates the need to work with dimensions greater than or equal to $3$. Consider the following arbitrary  $2 \times 2$ matrix
 
$${\bf X} = \left[\begin{matrix}
 a &b \cr c &d
\end{matrix} \right]\ .$$
If $a < b$ and $c < d$, then
$${\bf X}^{(1)} = \left[
\begin{matrix}
  -1 &1 \cr
  -1 &1
\end{matrix}
 \right],$$ 
 so $(S^{(1)}_j)^2 = 0$. 

If $a,b,c,d$ are, say, \emph{iid} and have symmetric distributions, then
$$P( (S^{(1)}_j)^2 = 0 ) = {1\over 2}.$$
 
A moment's reflection shows that if ${\bf X}$ is $I \times 2$,with $n$ odd, then after each row normalization,
each column has an odd number of entries, each entry being $+1$ or $-1$.  However, each row has exactly one $+1$ and one $-1$. Thus $(S^{(n)}_i - S^{(n+1)}_j)^2 \rightarrow 0$ is impossible. With $\min (I, J) \ge 3$, both tend to 1 a.s. as $n \nearrow \infty$. The fact that $\{X_{ij}\}$ iid implies in particular that ${\bf X}$ is row and column exchangeable.  Thus, if
 
$${\bf X} = [{\bf x}_{1, c}, ... {\bf x}_{J, c}],$$

\noindent where each column is $I \times 1$, and $\pi$ is a permutation of the integers $\{ 1, ..., J\}$,
then ${\bf X} \sim [{\bf x}_{\pi (1), c}, ..., {\bf x}_{\pi (J), c}]$.
The analogous holds for rows, where
 
$${\bf X} = \left[
\begin{matrix}
           {\bf x}_{1, r} \cr
                   \vdots \cr
                   {\bf x}_{I, r}
                 \end{matrix}
\right] ;$$ 

where ${\bf x}_{i,r}$ is $1 \times J$

\section{Background and Motivation}

A first step in the argument of \cite{olshen_rajaratnam1} is to note that Lebesgue measure on the Borel subsets of $\mathbb{R}^{IJ}$ is mutually absolutely continuous with respect to IJ product Gaussian measure, each coordinate being standard Gaussian.  Thereby, the distinction between measure and topology is blurred; arguments of \cite{olshen_rajaratnam1} as corrected here.  Having translated a problem concerning Lebesgue measure to one concerning Gaussian measure, one cannot help note from graphs in \cite{olshen_rajaratnam1} Figures 1, 2, 3, 6, and 7 that suggest with entries of $\bf{X}$ chosen independently from a common absolutely continuous distribution, and as led to these figures, almost surely ultimately, convergence is at least exponentially fast. In these graphs, the ordinate is always the difference of logarithms (the base does not matter) of squared Frobenius norm of the difference between current iteration and the immediately previous iteration; always abscissa indexes iteration. One purpose of this paper is to demonstrate the ultimately almost sure rapidity of convergence.  Readers will note we  assume that coordinates are independent and identically distributed, with standard normal distributions.  This assumption is used explicitly though unnecessarily in \cite{olshen_rajaratnam1}, but only implicitly here, and then only to the extent that arguments here depend on those in \cite{olshen_rajaratnam1}. Obviously, this Gaussian assumption is sufficient. It is pretty obviously not necessary.  

For the sake of motivation we illustrate the patterns of convergence from successive normalizations on $5\times5$ matrices. We first describe the details of the successive normalization in a manner very similar  to \cite{olshen_rajaratnam1} . Consider a simple $5\times5$ matrix with entries generated from a normal distribution with a given mean and standard deviation. In our case we take the mean to be 2 and variance to be 4, though the specific values the mean and variance parameter take do not really matter. We first standardize the initial matrix $\bf{X^{(0)}}$ at the level of each by row, i.e., first subtracting the row mean from each entry and then dividing each entry in a given row by its row standard deviation. The matrix is then standardized at the level of each column, i.e., by first subtracting the column mean from each entry and then by dividing each entry by the respective column standard deviation. Row mean and standard deviation polishing followed by column mean and standard deviation polishing is defined as one iteration in the process of attempting to row and column standardize the matrix. The resulting matrix is defined as $\bf{X^{(1)}}$. Now the same process is repeated with $\bf{X^{(1)}}$ and repeated until successive renormalization eventually yields a row and column standardized matrix. The successive normalizations are repeated until ``convergence" which for our purposes is defined as the difference in the squared Frobenius norm between two consecutive iterations being less than $10^{-8}$. 

The figures (see Figures 1 - 2) are plots of the log of the ratios of the squared Frobenius norms of the differences between consecutive iterates. In particular, they capture the type of convergence patterns that are observed in the $5\times5$ case from different starting values. In \cite{olshen_rajaratnam1} it was proved that regardless of the starting value (provided the dimensions are at least 3) the process of successive normalization always converges, in the sense that it leads to a doubly standardized matrix. In addition, it was noted that the convergence is very rapid. We can see empirically from Figures 1-11 that eventually the log of the successive squared differences tends to decay in a straight line, i.e., the rate of convergence is perhaps exponential. This phenomenon of linear decay between successive iterations in the log scale is observed in all the diagrams. Hence, one is led naturally to ask whether this is always true in theory, and if so under what conditions. The rate of convergence of successive normalization and related questions are addressed in this paper. 

\begin{figure}[!t]
\centering
\includegraphics[width=0.495\textwidth]{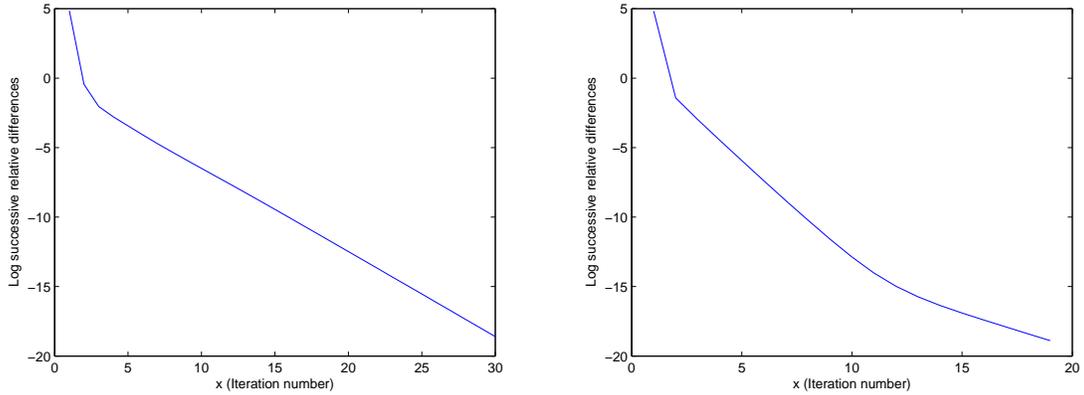}
\includegraphics[width=0.495\textwidth]{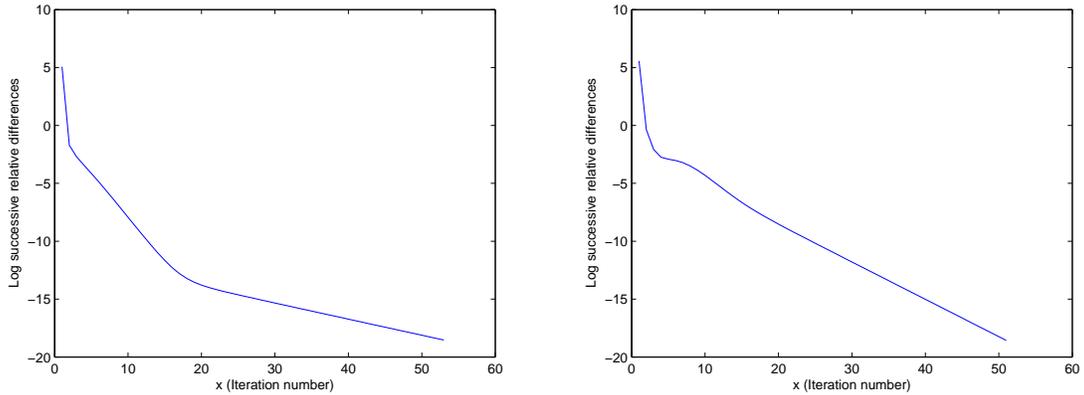}
\caption{Examples of patterns of convergence}
\label{fig_sim}
\end{figure}





\begin{figure}[!t]
\centering
\includegraphics[width=0.495\textwidth]{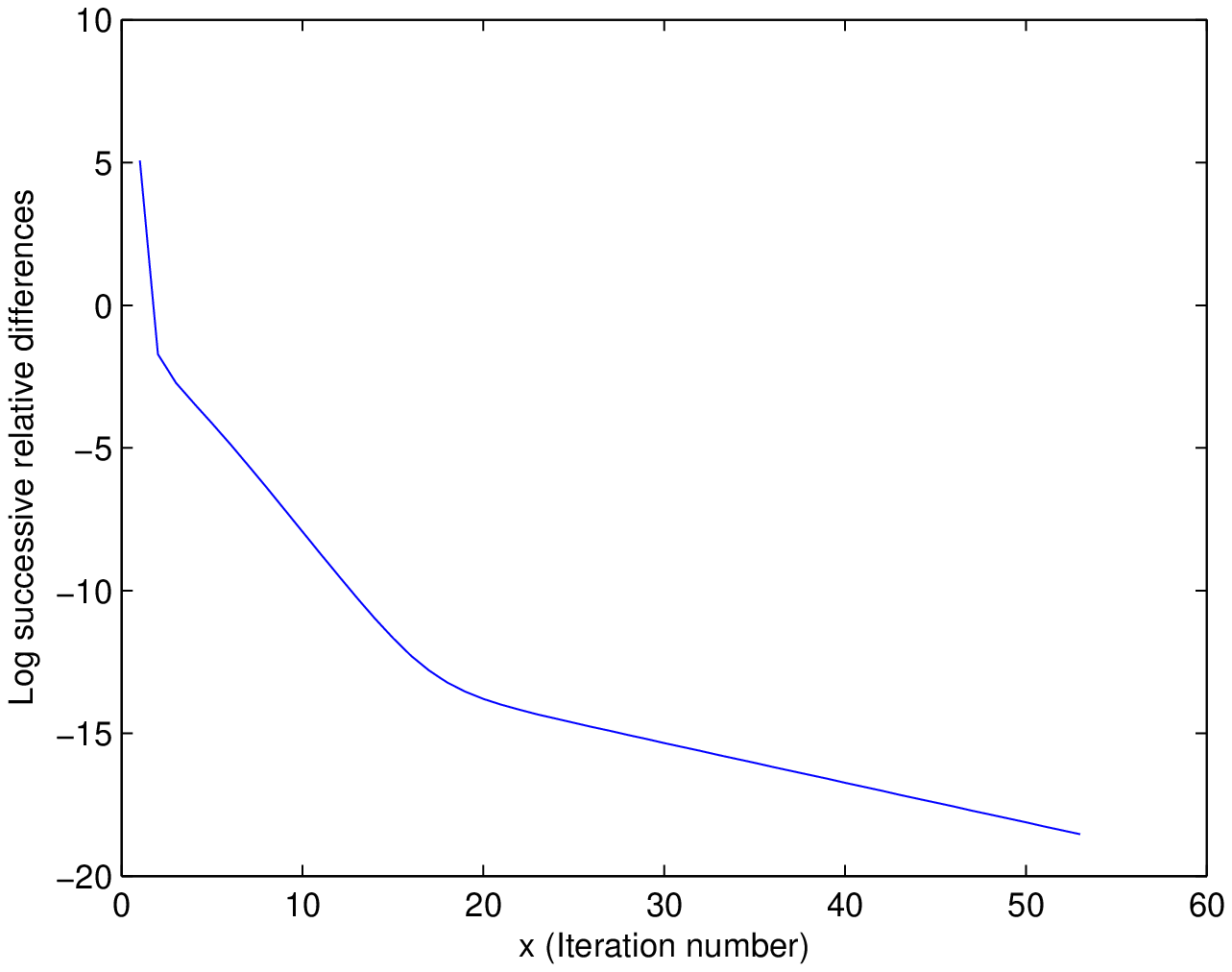}
\includegraphics[width=0.495\textwidth]{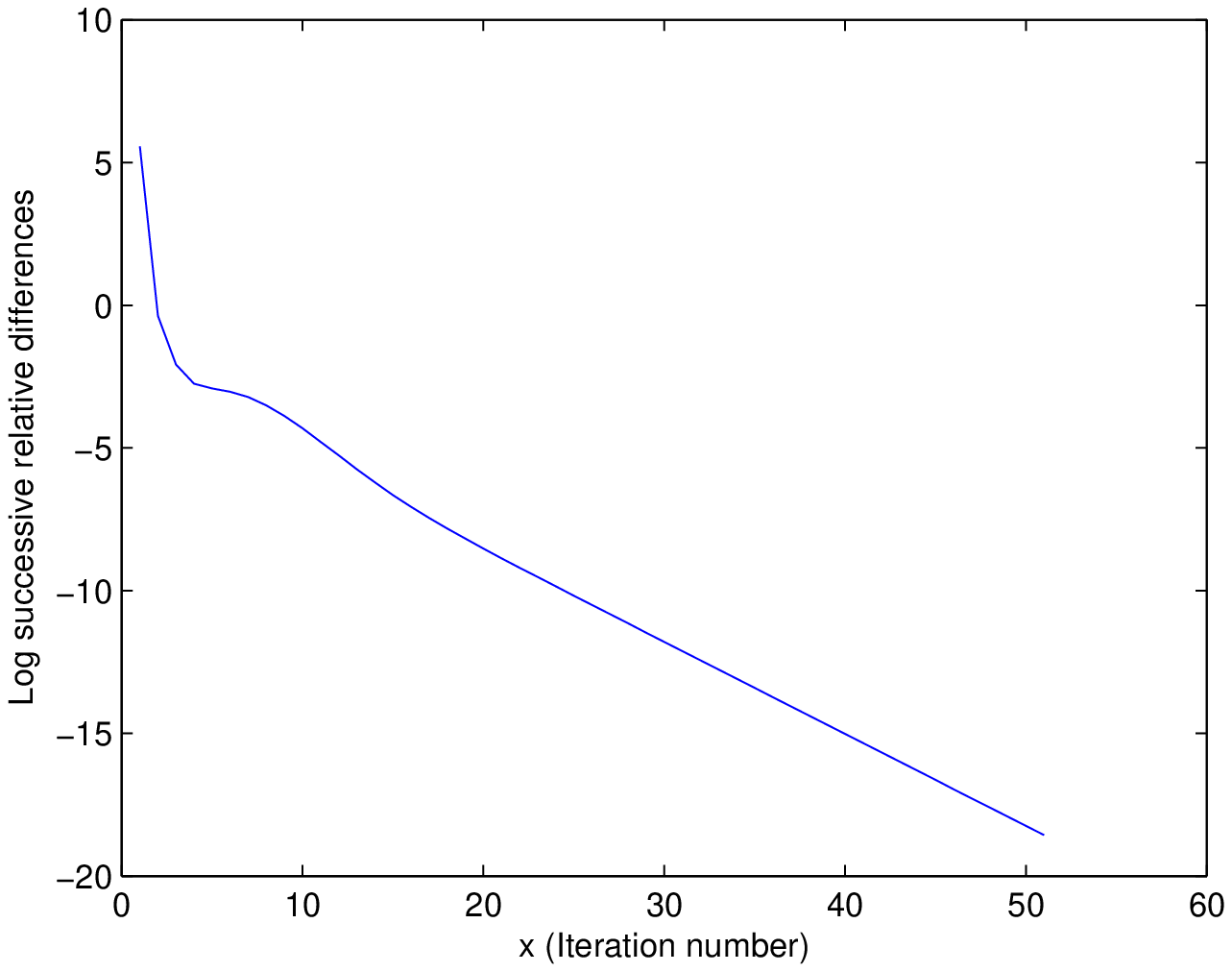}
\caption{Examples of patterns of convergence}
\label{fig_sim}
\end{figure}

A natural question to ask is whether the convergence phenomenon observed will still occur if simultaneous normalization is undertaken as compared to successive normalization. In other words the row and column mean polishing and row and column standard deviation polishing is all done at once: 

$$ X_{ij}^{t+1} = \frac{X_{ij}^{t} - {\bar X}_{i.}^{t} - {\bar X}_{.j}^{t}}{s_{i.} s_{.j}}$$

In this case the simultaneous normalization algorithm does not converge. In fact it is shown to diverge. Figures 3-8 below illustrate this phenomenon.

\begin{figure}[!t]
  \centering
  \includegraphics[width=0.495\textwidth]{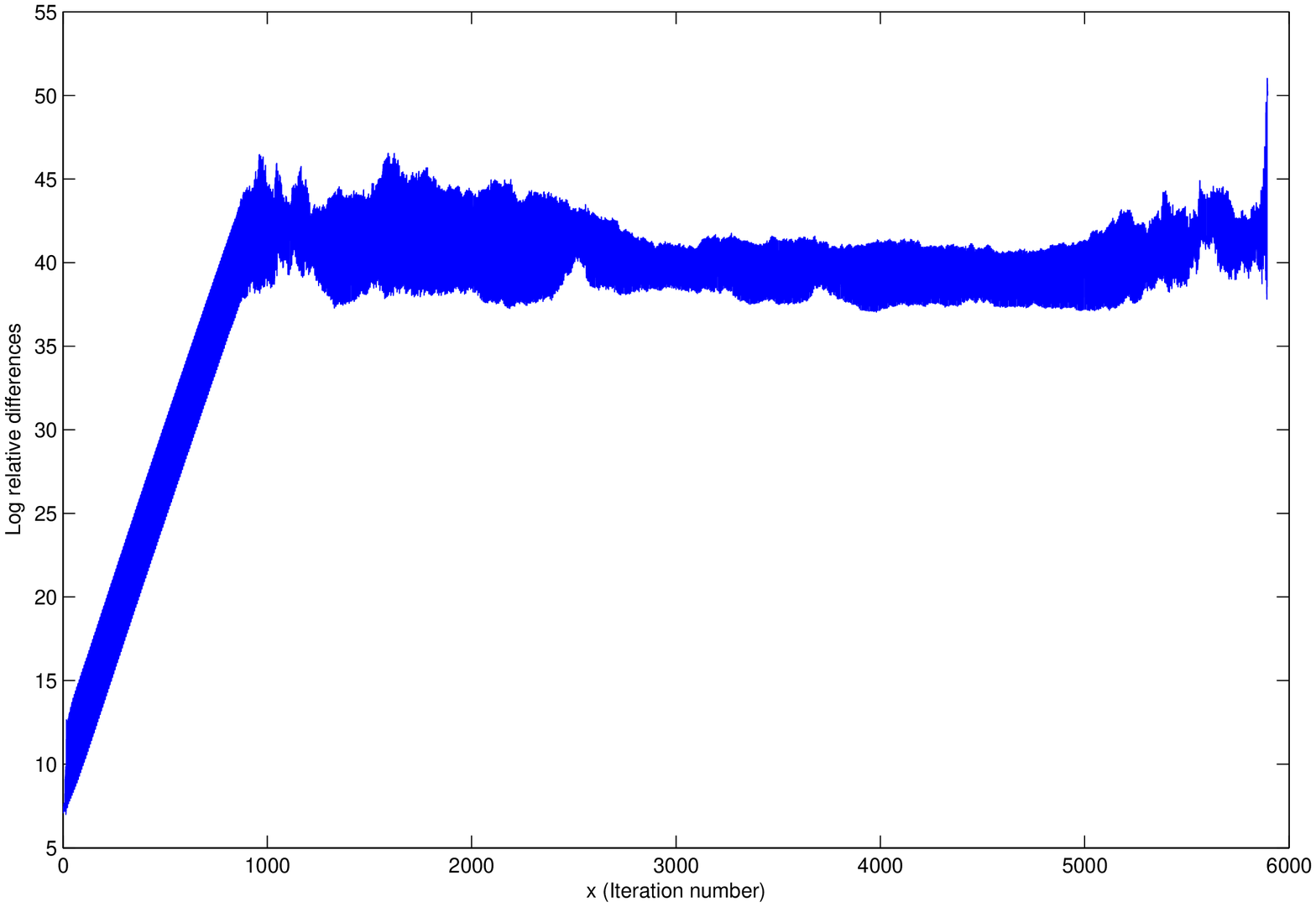}
  \includegraphics[width=0.495\textwidth]{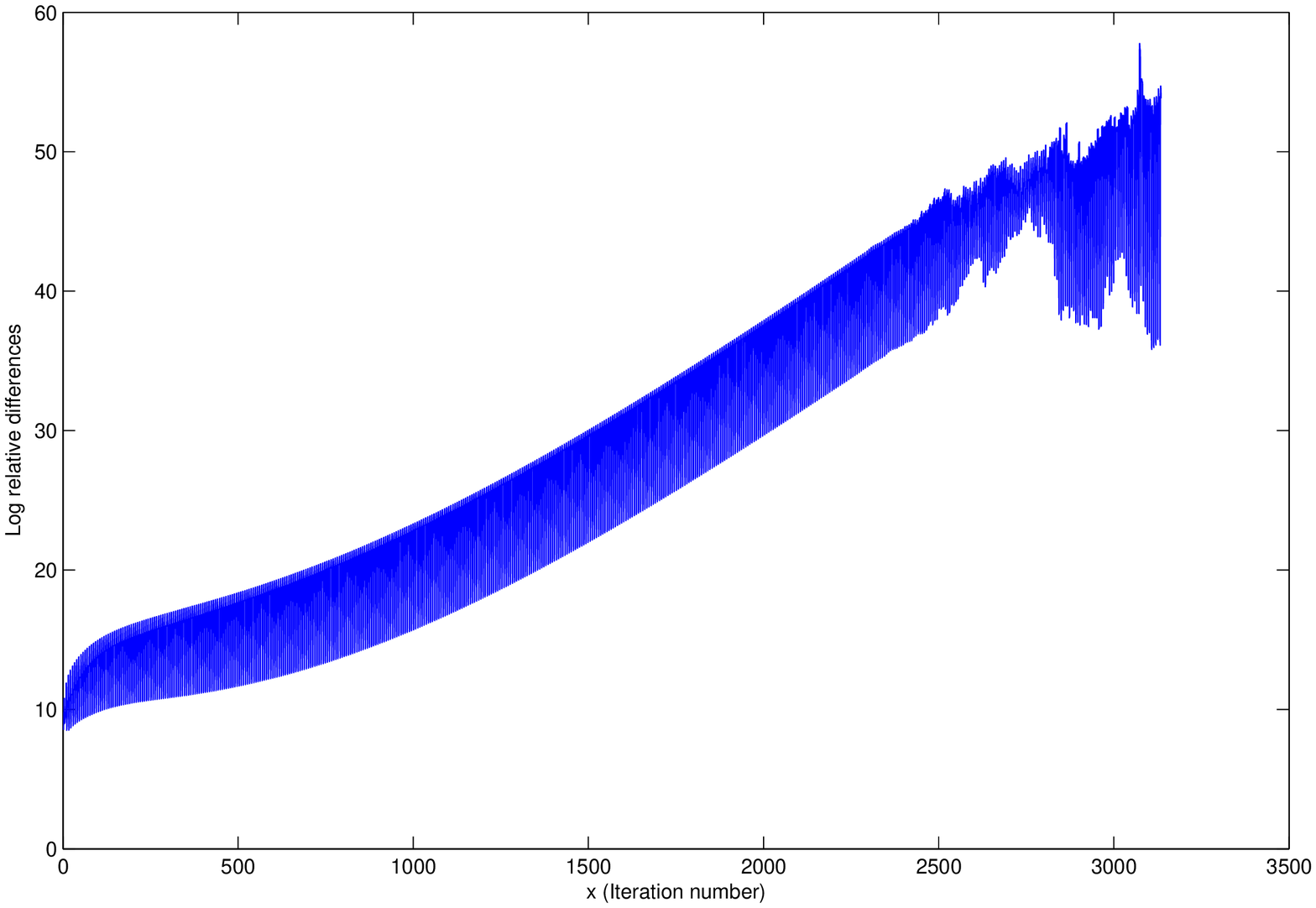}
  \caption{Two examples of simultaneous normalization: these illustrate that simultaneous normalization does not lead to convergence}
\end{figure}
\begin{figure}[!t]
  \centering
\includegraphics[width=0.495\textwidth]{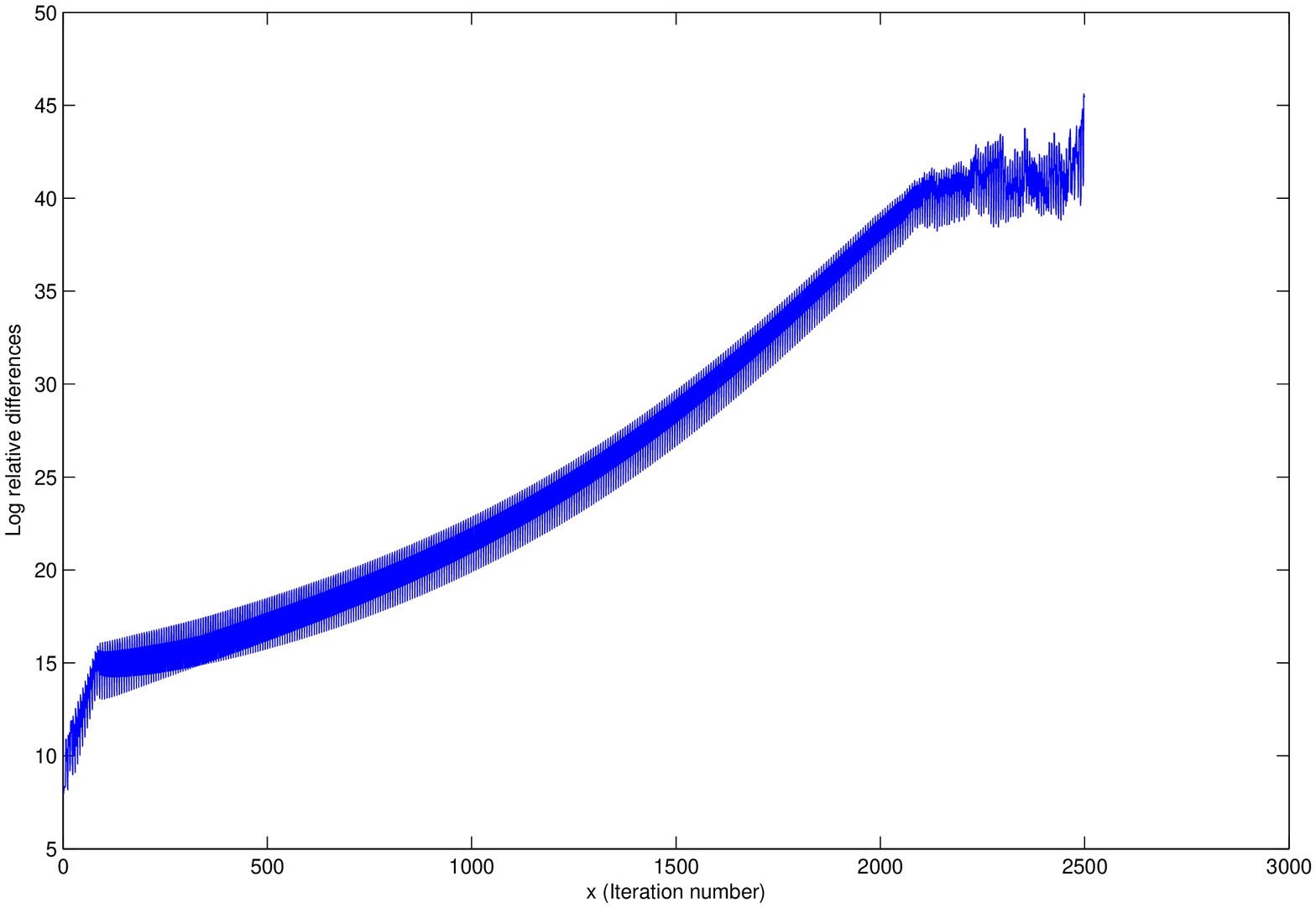}
\includegraphics[width=0.495\textwidth]{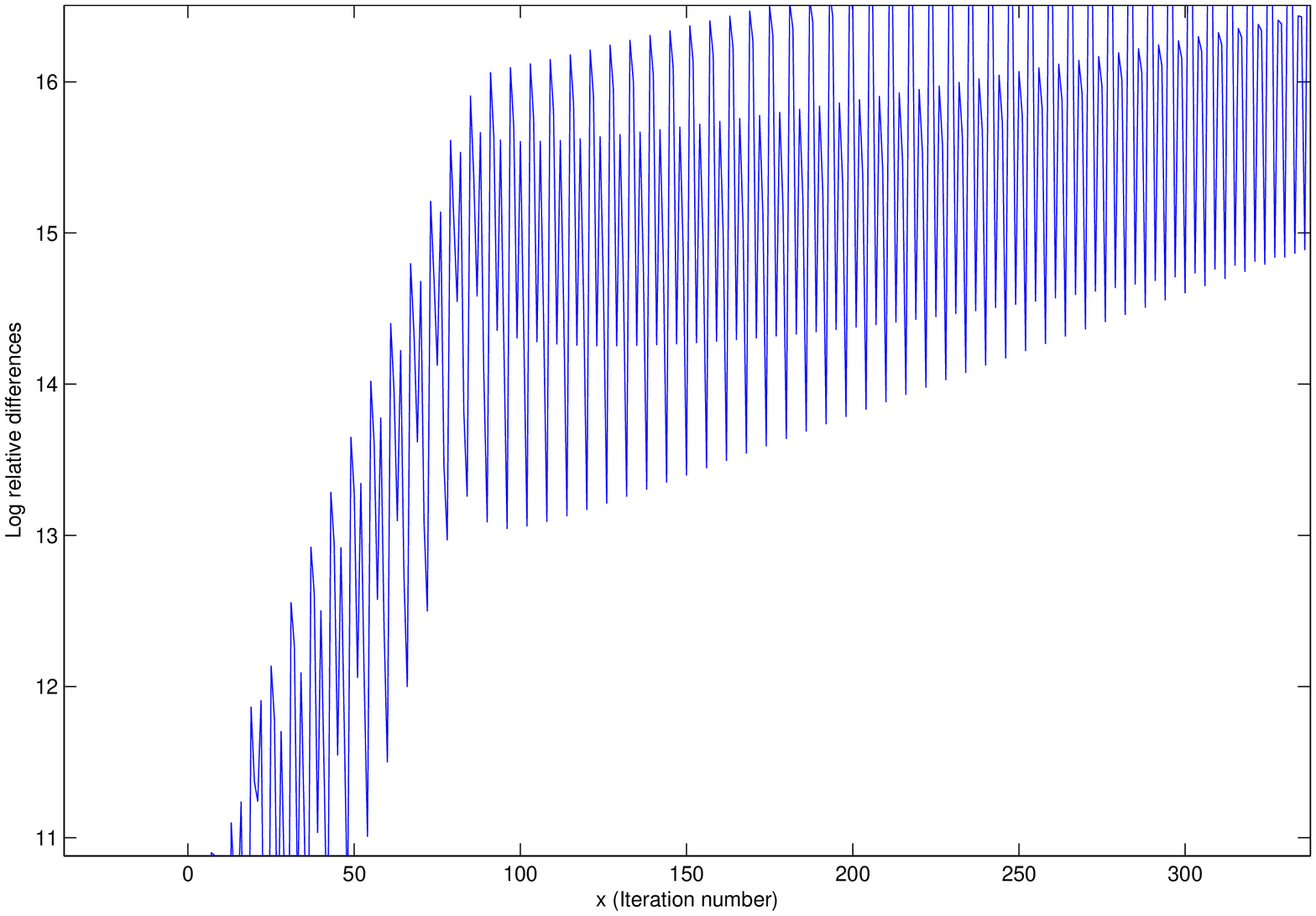}
  \caption{(Left) Example of simultaneous normalization, (Right) Same example zoomed in. }
\end{figure}
\begin{figure}[!t]
  \centering
\includegraphics[width=0.495\textwidth]{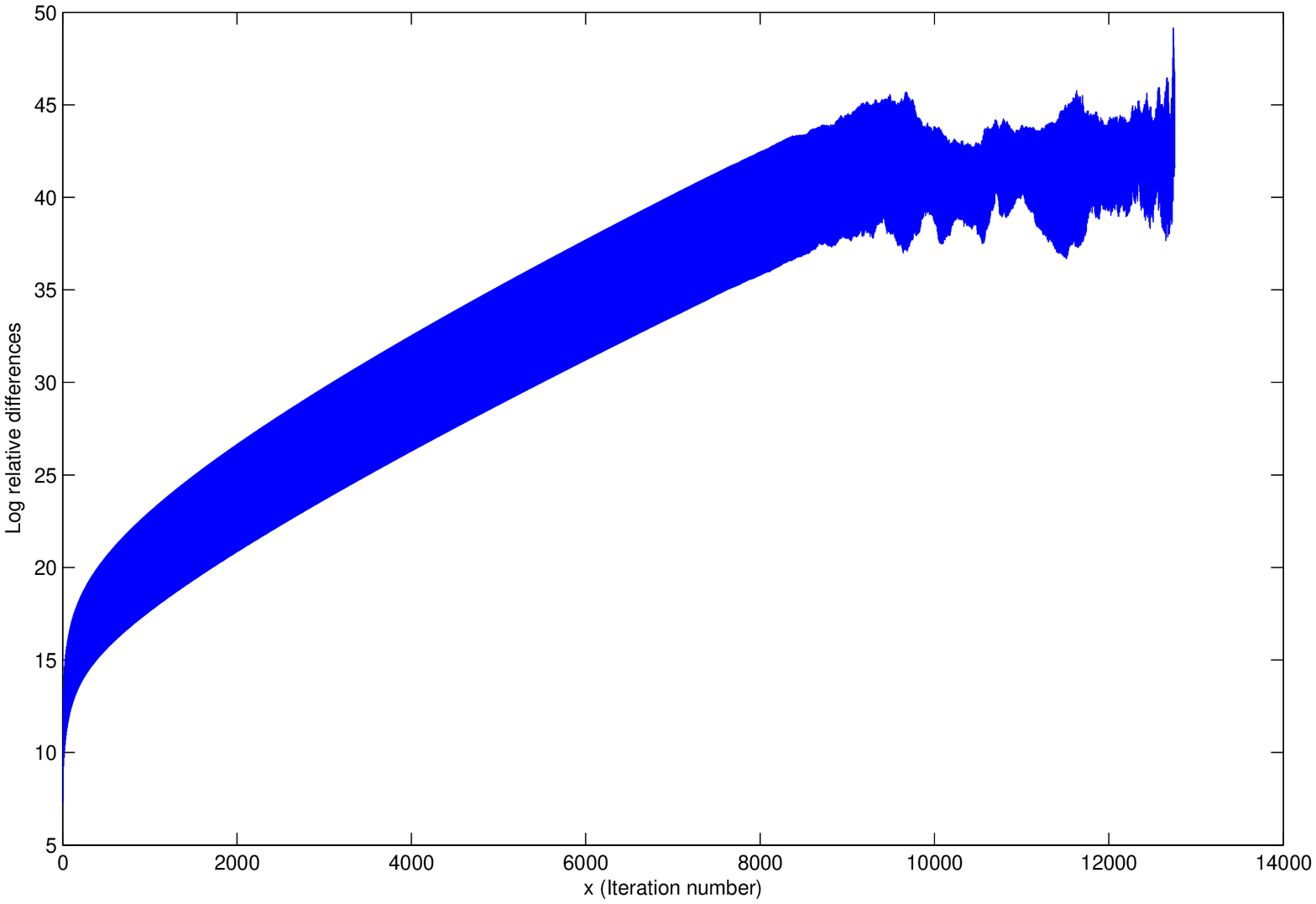}
\includegraphics[width=0.495\textwidth]{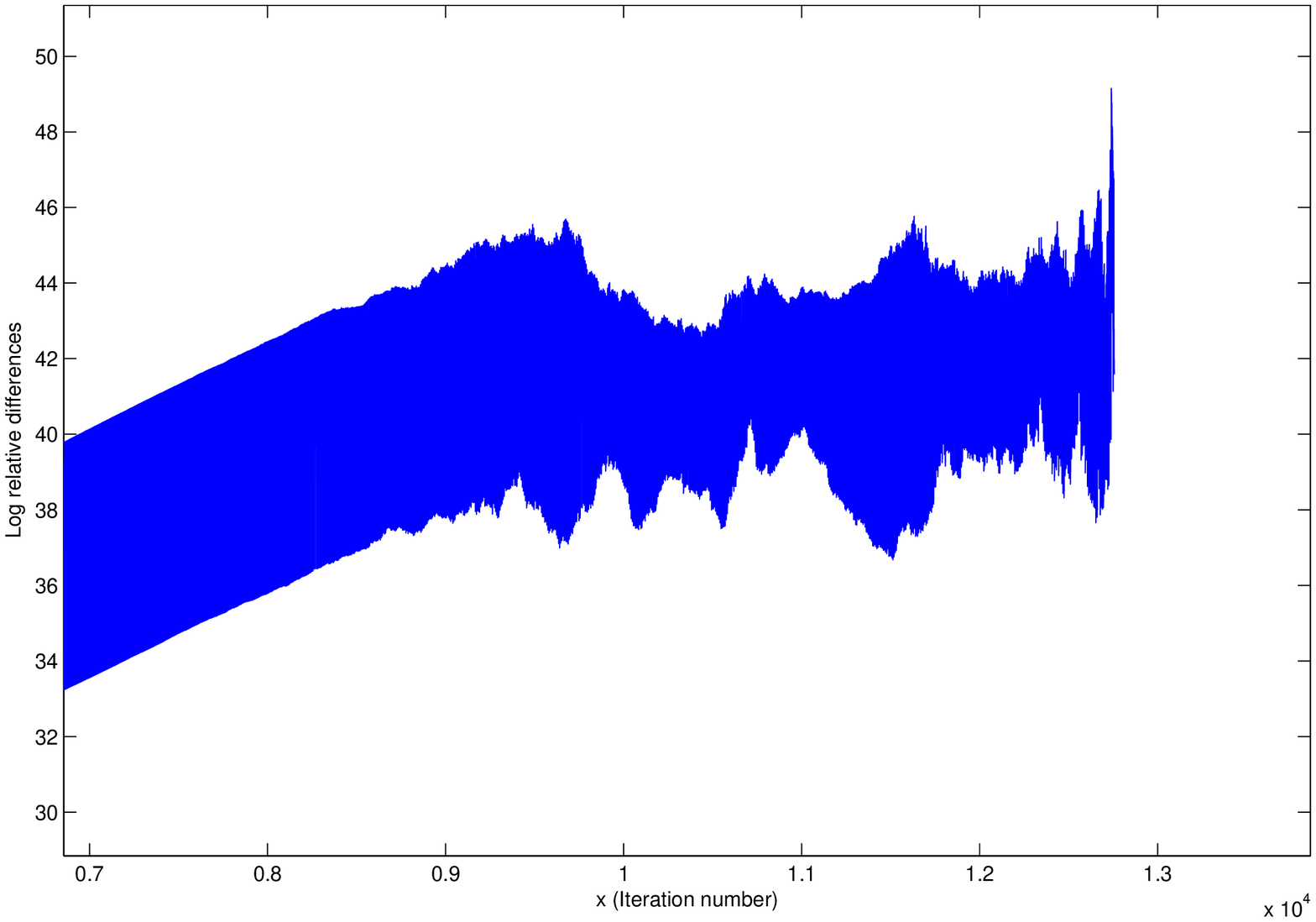}
  \caption{(Left) Example of simultaneous normalization, (Right) Same example zoomed in. }
\end{figure}
\begin{figure}[!t]
  \centering
\includegraphics[width=0.495\textwidth]{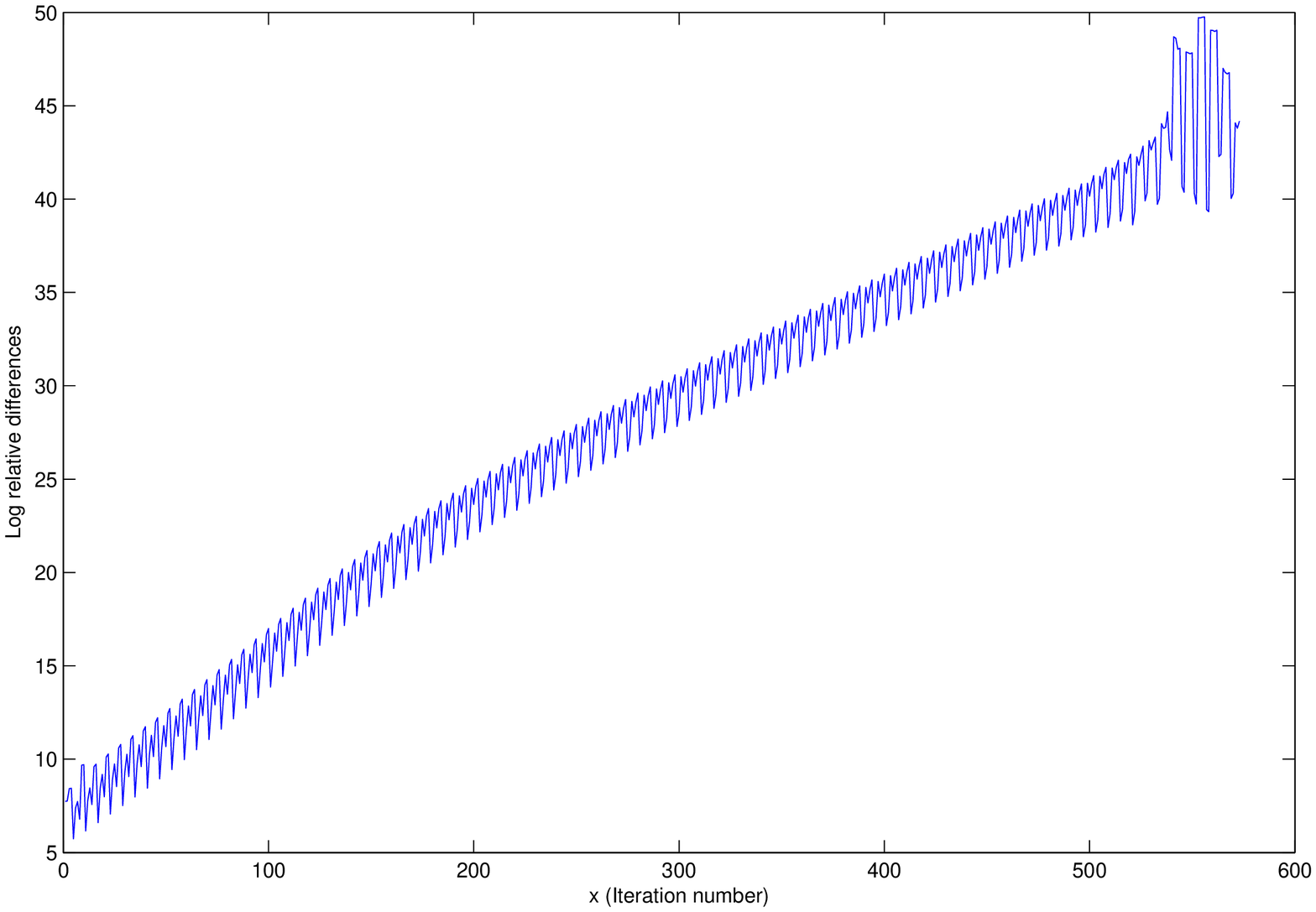}
\includegraphics[width=0.495\textwidth]{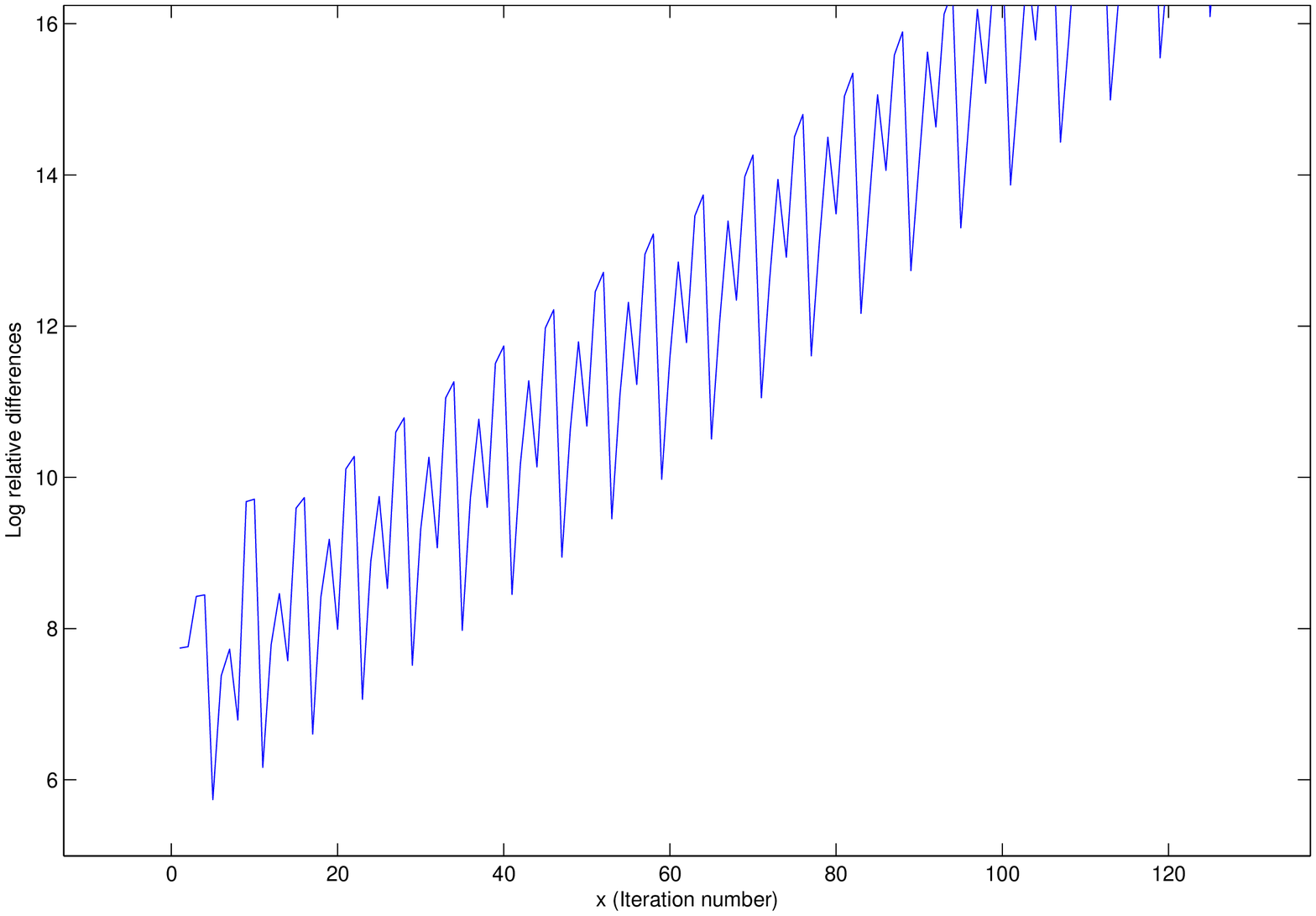}
  \caption{(Left) Example of simultaneous normalization, (Right) Same example zoomed in. }
\end{figure}
\begin{figure}[!t]
  \centering
  \includegraphics[width=0.495\textwidth]{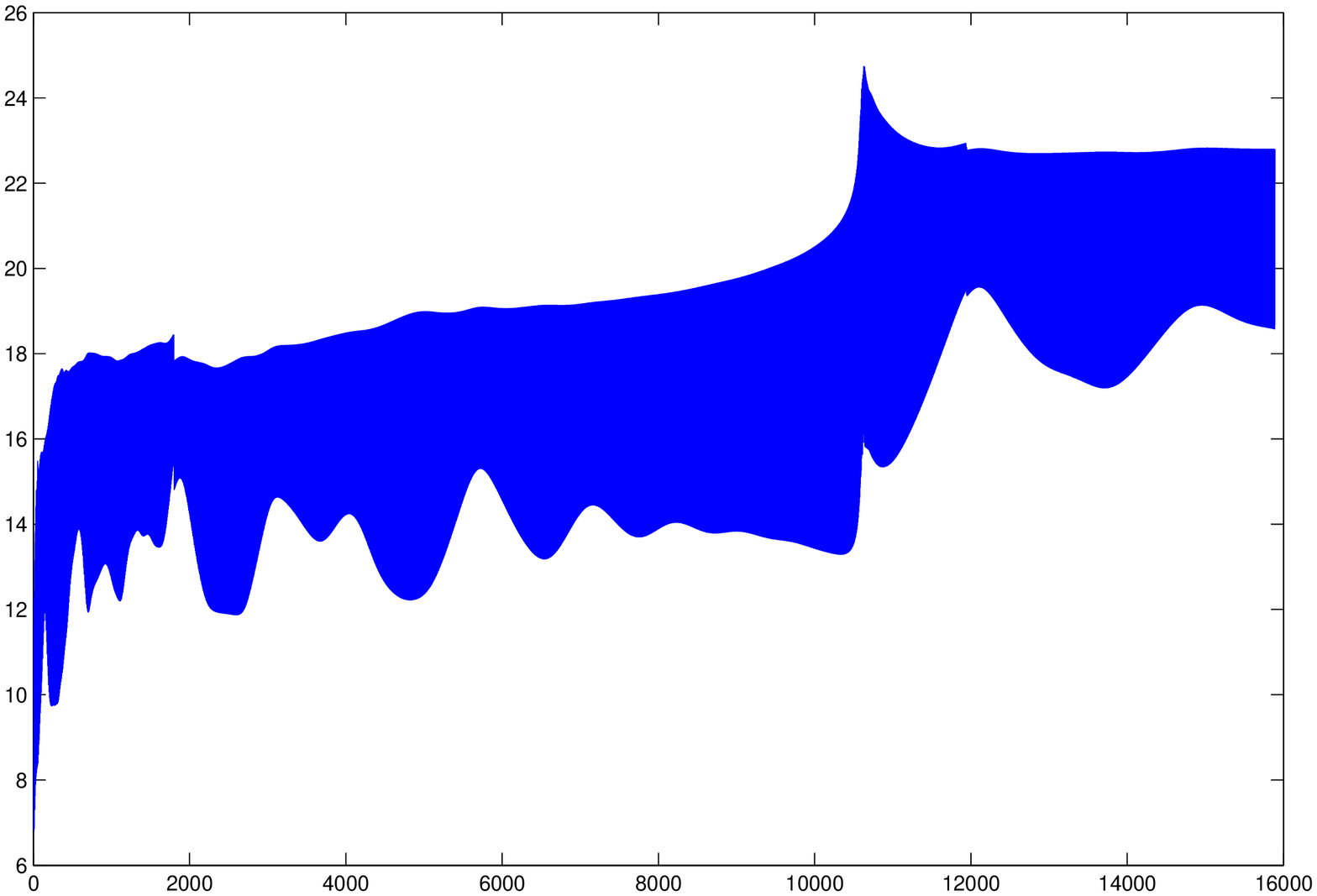}
  \includegraphics[width=0.495\textwidth]{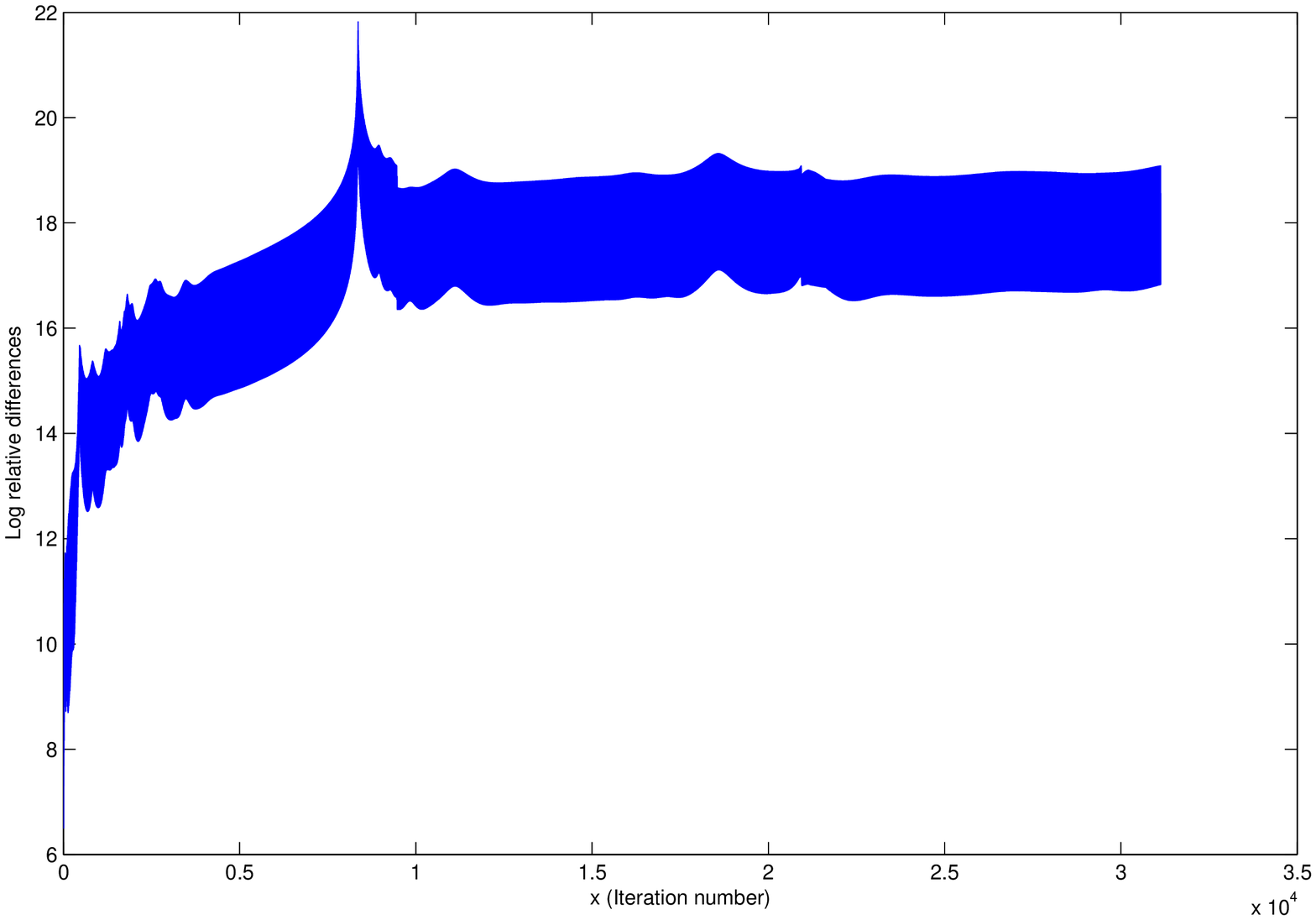}
  \caption{Two further examples of simultaneous normalization: these illustrate that simultaneous normalization does not lead to convergence}
\end{figure}
\begin{figure}[!t]
  \centering
  \includegraphics[width=0.495\textwidth]{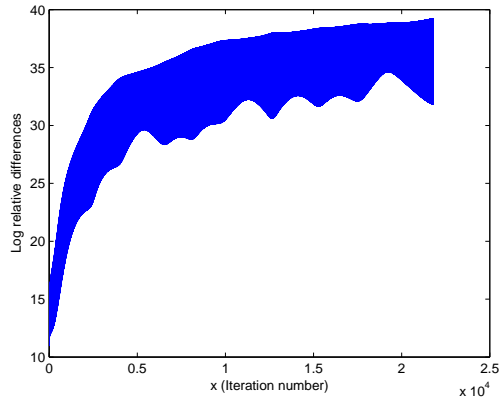}
  \caption{Example of simultaneous normalization}
\end{figure}

\newpage
\section{Convergence and Rates}

\noindent
Theorem 4.1 of \cite{olshen_rajaratnam1} is false.  A claimed backwards martingale is NOT.  Fortunately, all that seems damaged by the mistake is pride.  Much is true.  To establish what this is requires notation. \ \
\vskip.3cm
 
\noindent
$\Xb$ denotes an $I\times J$ matrix with values $x\in\Re^{IJ}$. We take coordinates $X_{ij}(x)=\Xb_{ij}$ to be iid $N(0,1)$. As in \cite{olshen_rajaratnam1}, \cite{olshen_rajaratnam2}, $3\leq \min(I,J)\leq\max(I,J)<\infty$. As before,
%
%
%
\begin{eqnarray}
  (S_i^{(0)})^2 &= & {1\over J} \sum_{j=1}^{J} (X_{ij}-\bar{X}_{i\cdot}^{(0)})^2 \\
  &= & {1\over J} \sum_{j=1}^{J} (X_{ij})^2 - {2\over J} \sum_{j=1}^{J} X_{ij}\bar{X}_{i\cdot}^{(0)} +
  (\bar{X}_{i\cdot}^{(0)})^2  \\
  &= & {1\over J}\sum_{j=1}^{J} (X_{ij})^2 - 2\bar{X}_{\cdot j}\bar{X}_{i\cdot} + (\bar{X}_{i\cdot})^2.
\end{eqnarray}
$\Xb^{(1)} = [X_{ij}^{(1)}]$, where, in an obvious notation $X_{ij}^{(1)} = (X_{ij}-\bar{X}_{i\cdot}^{(0)})/S_i^{(0)}$.
 
\vskip.3cm
\noindent
In view of \cite{olshen_rajaratnam1}, almost surely $(S_i^{(0)})^2 > 0$. By analogy, set $\Xb^{(2)}=[X_{ij}^{(2)}]$, where $X_{ij}^{(2)}=(X_{ij}^{(1)}-\bar{X}_{\cdot j}^{(1)})/S_j^{(1)}$, the definition of $(S_j^{(1)})^2$ being clear. As in \cite{olshen_rajaratnam1}, $(S_j^{(1)})^2>0$. In general, for $m$ odd, $X_{ij}^{(m)}=(X_{ij}^{(m-1)}-\bar{X}_{i\cdot}^{(m-1)})/S_i^{(m-1)}$. 
\vskip.3cm
 
\noindent
Likewise, for $m$ even $X_{ij}^{(m)}=(X_{ij}^{(m-1)}-\bar{X}_{\cdot j}^{(m-1)})/S_j^{(m-1)}$. Without loss, we take $(S_i^{(l)})^2$ and $(S_j^{(l)})^2$ to be positive for all $(i,j,l)$.

In various places we make implicit use of a theorem of Skorokhod, the Heine-Borel Theorem, and this obvious fact.  Suppose we are given a sequence of random variables ${\bf Y} = (Y_1, Y_2, ...)$ and a condition ${\cal C}_0$ that depends only on their finite-dimensional distributions.  If we wish to make a conclusion ${\cal C}_1$ concerning ${\bf Y}$, then it is enough to find one probability space that supports ${\bf Y}$ with given finite dimensional distributions for which ${\cal C}_0$ implies ${\cal C}_1$.  The theorem of Skorokhod mentioned(see pages 6-8 of \cite{billingsley}) is to the effect that if the underlying probability space is (for measure theoretic purposes) the real line with measures given by, say, distribution functions ${\bf F} = (F_1, F_2, ...)$; and if ${\bf F}$ is compact with respect to weak convergence, where $F_{q_l}$ converges to $G$ in distribution; then if each $F_i$ is absolutely continuous with respect to some $H$ (which itself is absolutely continuous), and $Y_i(x)= F_{i}(x)$, then $Y_{q_l}$ converges H-almost surely to a random variable $Y$.
\vskip.5cm
 
\noindent
Off of a set of probability $0$, 
$\sum_{i,j} (X^{(q)}_{ij})^2 = I J$ for all $q \ge 1$.  We assume that ${\bf x}$ lies outside this ``cursed'' set. This is key to convergence.\\
 
\noindent
Almost surely, $(S^{(2m-1)}_j)^2$ has positive limsup as $m$ increases without bound.\\
 
\noindent
Let $A_j = \{ \overline{\lim}_m (S_j^{(2m-1)})^2 = 0\}$. \\
 
\noindent $P(A_j) = 0$ \ \ \ $j = 1,2, ... \ J$

Since the entries of ${\bf X}$ are independent, ${\bf X}$ is row and column exchangeable.  This property is inherited by ${\bf X}^{(q)}$ for every $q$.  Because all entries (for $q \ge 1$) are a.e. bounded uniformly in ${\it q}$, $E\{ X_{ij}^{(q)}\}$ and $E\{ (X_{ij}^{(q)})^2 \}$ exist and are finite (with fixed bound that applies to all $q$).  Exchangeability implies that all $E\{ X_{ij}^{(q)} \} = 0$, and all $E\{ (X_{ij}^{(q)})^2 \} = 1$.  Bounded convergence implies that if $(S_j^{(2m-1)})^2$ tends to $0$ along a subsequence as $m$ increases, not only is the limit bounded as a function of ${\bf x}$, but also the limit random variable has expectation $0$.  Necessarily every almost sure subsequential limit in $m$ of the random variables $\overline{X}_{\cdot j}^{(2m-1)}$ has mean $0$.  Likewise, every almost sure subsequential limit in $m$ of the random variables $(X_{ij}^{(2m-1)})^2$ has expectation 1.  All are bounded as functions of ${\bf x}$.  One consequence of these matters is that $P(A_j)=0$. The next paragraph is a proof.
\vskip.5cm
 
\noindent
The only subsequential almost sure limits of
$\{(S_j^{(2m-1)})^2: j=1, 2, ... \}$ and of
$\{(S_j^{(2m)})^2: j=1, 2, ... \}$ have expectation 1.  Fix a $j, 1 \le j \le J$.  Let $E=E(j) = \{i:$ for some $q_l, P(\overline{\lim} | X_{ij}^{(q_l)} | > 0) > 0\}$.  Column exchangeability implies that for any pair $(i, i'), i\in E(j)$ iff $i' \in E(j)$.  The first sentence entails that $E(j)$ is not empty.  Therefore, $E(j) = \{ 1, ..., I\}$. \ \
\vskip.4cm
 
\noindent
Let $i_0 \not= i_1$.  There is a subsequence of $\{ q_l\}$ -- for simplicity write it as $\{ q_l\}$ -- along which almost surely
$$(a) \; \lim_{q_l} X_{i_o j}^{(2_{q_l} -1)}  \ \mbox{and} \  \lim_{q_l} X_{i_0 j}^{(2_{q_l})}\ \mbox{both exist;} $$
$$(b) \; \lim_{q_l} X_{i_1 j}^{(2_{q_l} -1)}  \  \mbox{and} \  \lim_{q_l} X_{i_1 j}^{(2_{q_l})}\ \mbox{both exist; and} $$

\noindent If $ P(\lim_{q_l} | X_{i_o j}^{(2_{q_l})} - X_{i_1 j}^{(2_{q_l})} | = 0) = 1 $ then exchangeability implies that $(S_j^{(2_{q_l} -1)})^2 \rightarrow 0$. So without loss of generality, if $E = \{\overline{\lim_{q_l}} | X_{i_o j}^{(2_{q_l})} - X_{i_1 j}^{(2_{q_l})} | > 0 \}$ then $P(E)> 0$. Write

$$(X_{i_o j}^{(2_{q_l})} - X_{i_o j}^{(2_{q_l} -1)}) - (X_{i_1 j}^{(2_{q_l})} - X_{i_1 j}^{(2_{q_l} -1)}) =
 (X_{i_o j}^{(2_{q_l} -1)} - X_{i_1 j}^{(2_{q_l} -1)}) (S_j^{(2_{q_l} -1)} - 1) / S_j^{(2_{q_l} -1)}.$$ 

Since $(X_{i_o j}^{(2_{q_l})})^2 - (X_{i_o j}^{(2_{q_l} -1)})^2 \rightarrow 0$ a.s.((a) and the expectations of both are 1), and likewise with $i_o$ replaced by $i_1$, ((b), etc..), and because $x^2 - y^2 = (x-y)(x+y)$, the first expression of the immediately previous display tends to $0$ on $E$.  So, too, does the second expression.  This is possible only if $S_j^{(2_{q_l} -1)} \rightarrow 1$ on $E$ (we take the positive square root).  On $E^c$, $S_j^{(2_{q_l} -1)} \rightarrow 0$. Since $E[S_j^{(2_{q_l} -1)}] \rightarrow 1$, $P(E^c)=0$. Further,
%
%
%
%
%

$$(X_{i_o j}^{(2_{q_l})} - X_{i_o j}^{(2_{q_l} -1)})  =
  { {X_{i_o j}^{(2_{q_l} -1)} (S_j^{(2_{q_l} -1)} - 1) + \overline{X}_{\cdot j}^{(2_{q_l} -1)} }
  \over {S_j^{(2_{q_l} -1)}}
  }
$$
\vskip.5cm
 
\noindent
As a corollary, one sees now that $\overline{X}_{\cdot j}^{(2_{q_l} - 1)} \rightarrow 0$ a.s.  Since the original $\{q_l\}$ could be taken to be an arbitrary subsequence of $\{q\}$, we conclude that convergence of row and column means to $0$ and convergence of row and column standard deviations to 1 takes place everywhere except on a set of Lebesgue measure $0$.
\vskip.4cm
 
\noindent          
THEOREM 4.1  Efron's algorithm converges almost surely for ${\bf X}$ on a Borel set of entries with complement a set of Lebesgue measure $0$.

\vspace{1cm}

\noindent We turn now to a study of rates of convergence. To begin, define the following

$${\bf Z} = \lim_m {\bf X}^{(m)} \; \textrm{ a.e. and } \lambda_{ij}^{(2m-1)} = X_{ij}^{(2m-1)} - Z_{ij}.$$ 

Now let  $\lambda_j = \max_i | \lambda_{ij}^{(2m-1)} |.$ Almost everywhere convergence and the fact that for each row and each column, not every $Z_{ij}$ can be of the same sign enable us to conclude that for $m$ large enough 

$$|\overline{X}_{\cdot j}^{(2m-1)} | \le { (I -1) \over I} \lambda_j^{(2m-1)}  $$
 
Remember that for $j = 1, ..., J \ge 3$, ${1\over I} \sum_{i=1}^I (Z_{ij})^2 = 1$, and analogously for $i = 1, ..., I \ge 3$. 

Now write
\begin{eqnarray*}
{1\over I} \sum_{i=1}^I (X_{ij}^{(2m-1)} - \overline{X}_{\cdot j} )^2 - {1\over I} \sum_{i=1}^I (Z_{ij} - 0)^2 
&=& a^2 - b^2 \\
&=& (a+b)(a-b),\;  \textrm{so} \\
a-1 &=& [ (a^2 - b^2) / a+ b] - 1  
\end{eqnarray*}
We know that For all $(i,j), \overline{Z}_{i \cdot} = \overline{Z}_{\cdot j} = 0$ a.e. To continue, we compute that

$$| X_{ij}^{(2m)} - X_{ij}^{(2m-1)} | =
      | X_{ij}^{(2m-1)} (S_j^{(2m-1)} -1) + \overline{X}_{\cdot j}^{(2m-1)}| \ / \ S_j^{(2m-1)},$$

where $S_j^{(2m-1)}$ is the positive square root of $(S_j^{(sm-1)})^2$.
 
Now write
\begin{eqnarray*}
(X_{ij}^{(2m-1)} - \overline{X}_{\cdot j}^{(2m-1)})^2 
&=& [(X_{ij}^{(2m-1)} - Z_{ij}) + (Z_{ij} - \overline{Z}_{\cdot j}) + (\overline{Z}_{\cdot j}- \overline{X}_{\cdot j})]^2\\
&=& [(X_{ij}^{(2m-1)} - Z_{ij})^2 + (Z_{ij} - \overline{Z}_{\cdot j})^2 + (\overline{Z}_{\cdot j}- \overline{X}^{(2m-1)})^2 \\
&+& \ 2 (X_{ij}^{(2m-1)} - Z_{ij}) (Z_{ij} - \overline{Z}_{\cdot j}) \\
&+& \ 2 (X_{ij}^{(2m-1)} - Z_{ij}) (\overline{Z}_{\cdot j} - \overline{X}_{\cdot j}^{(2m-1)}) \\
&+& \ 2 (Z_{ij} - \overline{Z}_{\cdot j}) (\overline{Z}_{\cdot j} - \overline{X}_{\cdot j}^{(2m-1)}) 
\end{eqnarray*}

\begin{figure}[!t]
  \centering
  \includegraphics[width=0.5\textwidth]{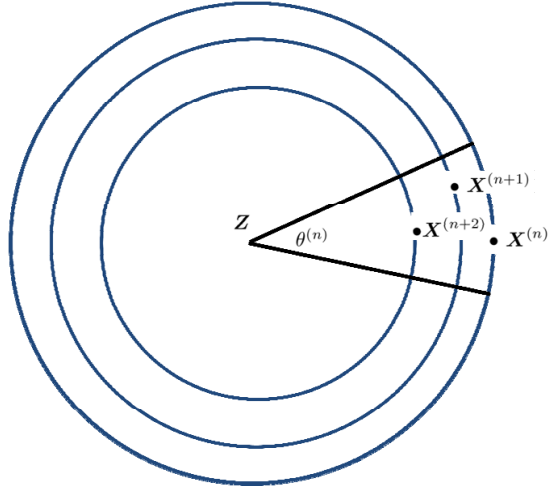}
  \caption{For $n$ large enough $\log\{{\|{\bf{X}}^{(n+2)} - {\bf{X}}^{(n+1)}\|}_F^2/{{\|\bf{X}}^{(n+1)} - {\bf{X}}^{(n)}\|}_F^2\} < 1 $ uniformly }
  \label{concentric_circles}
\end{figure}

One argues that
$${1\over I} \sum_{i=1}^I (X_{ij}^{(2m-1)} - \overline{X}_{\cdot j}^{(2m-1)})^2 \le D \lambda_j^2 + {2\over I} \sum_{(i=1)}^I (X_{ij}^{(2m-1)} - Z_{ij}) (Z_{ij} - \overline{Z}_{\cdot j}),$$

where $D = D(I, J) < \infty.$
A key observation is that for every $m$, there exists $j=j(m)$ for which $\{ \lambda_{ij}^{(2m-1)} \}$ are not of the same or all of the opposite sign as $\{ Z_{ij}: i=1, ..., I\}$, which, as was noted, are not, themselves of the same sign. Argue analogously regarding $\{ \lambda_{ij}^{(2m)}\}$. Refer now to Figure \ref{concentric_circles}. Our arguments show the correctness of the concentric circles in  $\mathbb{R}^{IJ}$ on which $\bf{X}^{(n)}, \bf{X}^{(n+1)}, \bf{X}^{(n+2)},\cdots $ ultimately lie. We conclude that to suitable approximation, the successive iterates lie on such circles with radii in geometric ratio. That ratio is uniformly $< 1$, but is \emph{not} arbitrarily close to $0$. However Figure \ref{concentric_circles} buries a key idea in our study of convergence are rates.  Thus, write 
$\theta^{(n)}$ for the angle between $X^{(n)}$ and $Z$. Obviously, the squared Frobenius norm $\|{\bf{X}}^{(n)} - {\bf{Z}}\|_F^2$ can be expressed as 

$$ \|{\bf{X}}^{(n)} - {\bf{Z}}\|_F^2 = \|{\bf{X}}^{(n)}\|_F^2 + \|{\bf{Z}}\|_F^2 - 
(2\cos(\theta^{(n)})\|{\bf{X}}^{(n)}\|_F \|{\bf{Z}}\|_F.$$

Now rewrite $\cos(\theta^{(n)}) = 1-\frac{1}{2\,IJ}(\|{\bf{X}}^{(n)} - {\bf{Z}}\|_F^2)$. Therefore, convergence of ${\bf{X}}^{(n)}$ to ${\bf{Z}}$ implies that $\theta^{(n)}$ can be taken arbitrarily close to $0$. This is another way of saying that for each $n$, $\|{\bf{X}}^{(n)}\|_F^2 = IJ$. Failure of this condition is why ``simultaneous normalization" fails.

\section{Illustrations of Convergence}

We now illustrate the rapidity of convergence in the $3\times3$ case to give a geometric perspective. First note that in the $3\times3$ case the set of fixed points are characterized by 3 unique values. For instance in the following doubly normalized matrix which arises a result of successive normalization has only three unique elements. 
\begin{eqnarray}
\bf{X}^{(\mathit{final})} &=&
\left[\begin{matrix}
-1.4137 & 0.7407  & 0.6730 \cr
 0.7407 & 0.6730 & -1.4137 \cr
0.6730 & -1.4137 & 0.7407
\end{matrix}\right]
\end{eqnarray}        
Hence we can use the first column of the limit matrix to represent the fixed point arising from successive normalization. Hence the curve of fixed points can be generated by applying the successive normalization process to random starting values. Figures 10 - 11 below illustrates the curve that characterizes the set of fixed points for the $3\times3$ matrix case when this numerical exercise is implemented. The origin is marked in black. The latter 3 subfigures superimpose the unit circle on the diagram in order to illustrate that the set of fixed points represent a ``ring" around the unit sphere. Figure 12 considers different starting values and their respective paths of convergence to the ``ring" of fixed points. Figures 13 - 15, presents 3 individual illustrations of different starting values and their respective paths of convergence to the ``ring" of fixed points. For each example, a magnification or close up of the path is provided. It is clear from these diagrams that the algorithm ``accelerates" or ``speeds up" as the sequence nears its limit.
\begin{figure}[!t]
\centering
\includegraphics[width=0.495\textwidth]{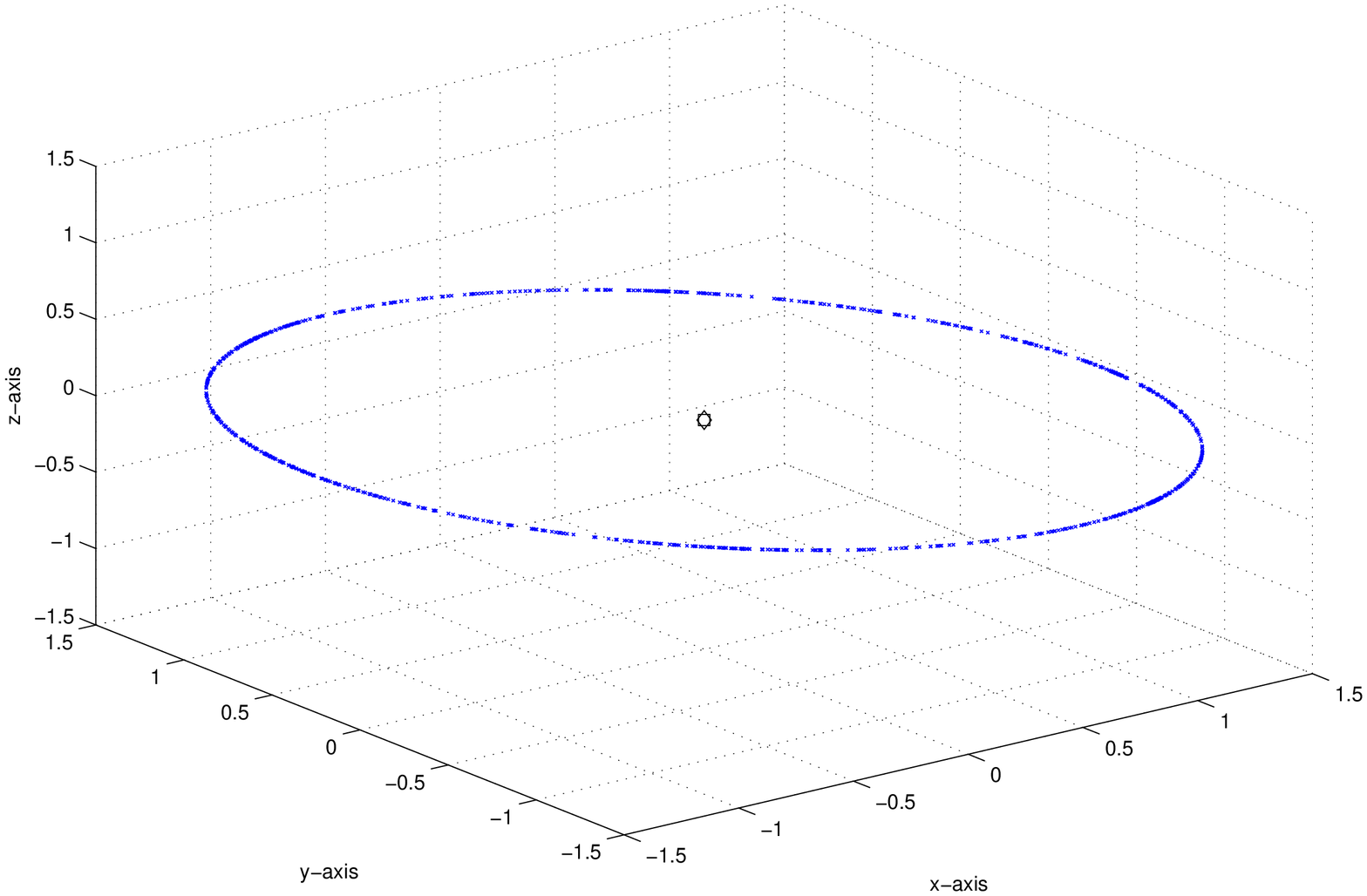}
\includegraphics[width=0.495\textwidth]{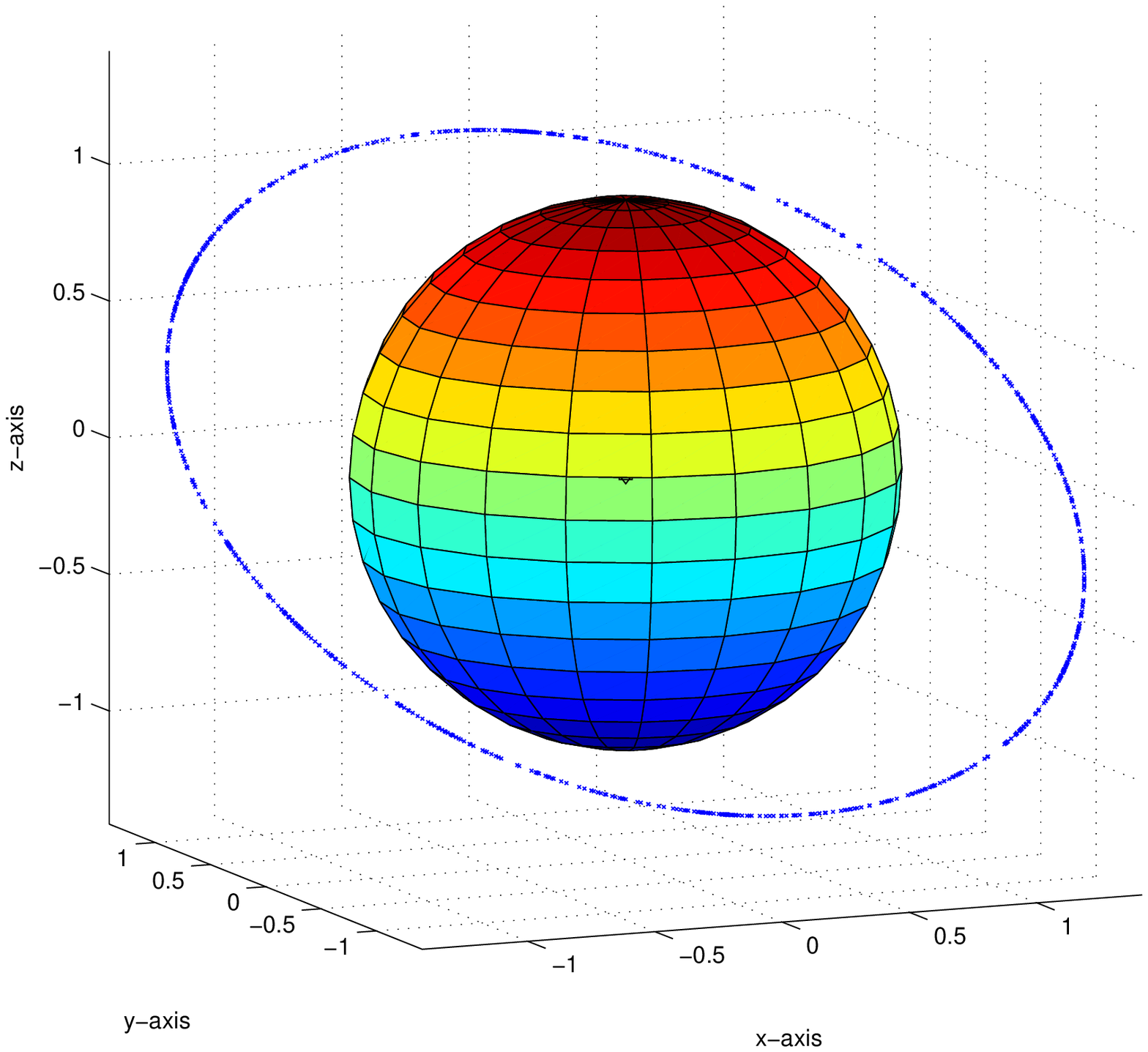}
\caption{(Left) Set of fixed points, (Right) Set of fixed points +
  unit sphere}
\label{fig_sim}
\end{figure}


\begin{figure}[!t]
\centering
\includegraphics[width=0.495\textwidth]{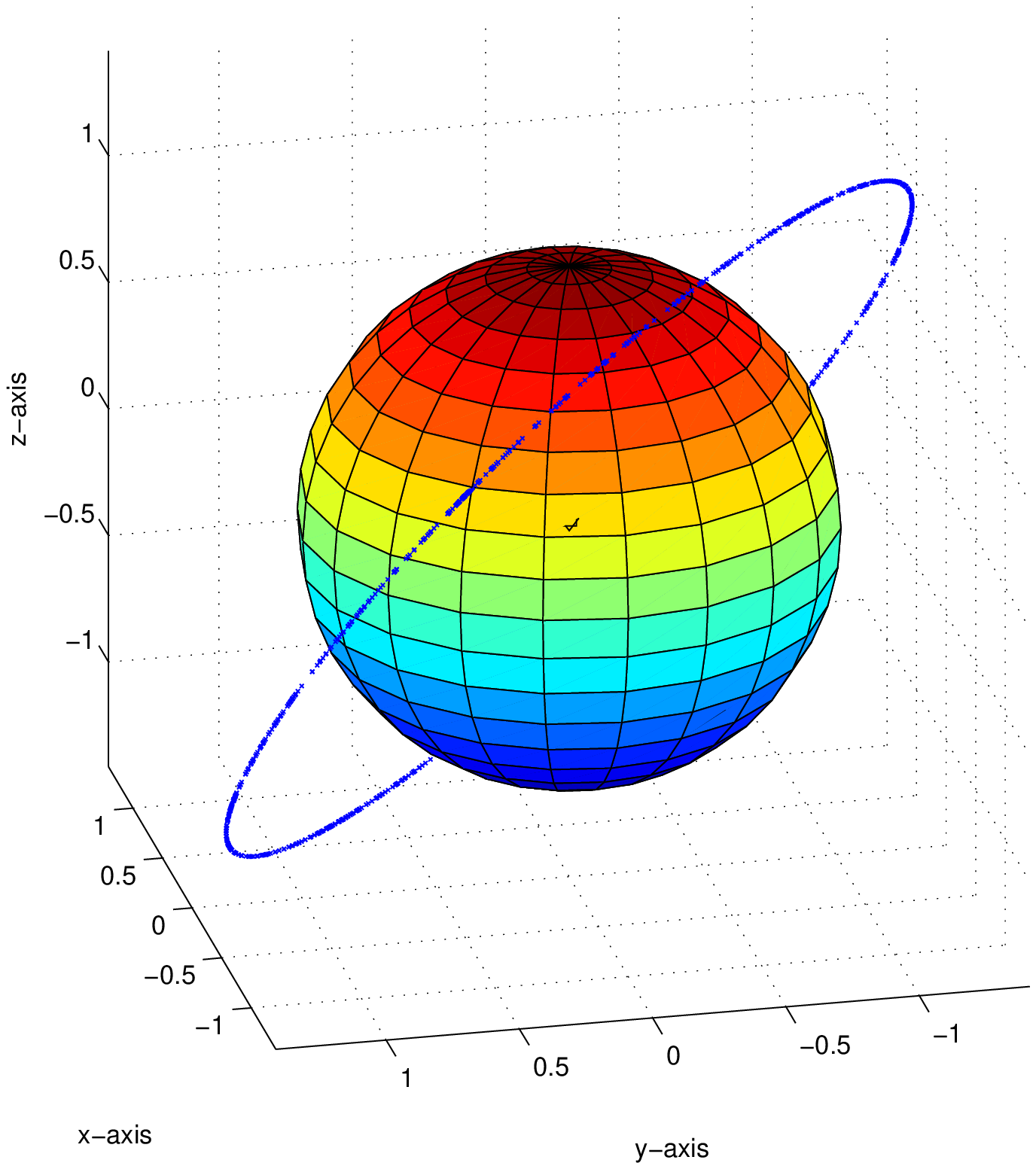}
\includegraphics[width=0.495\textwidth]{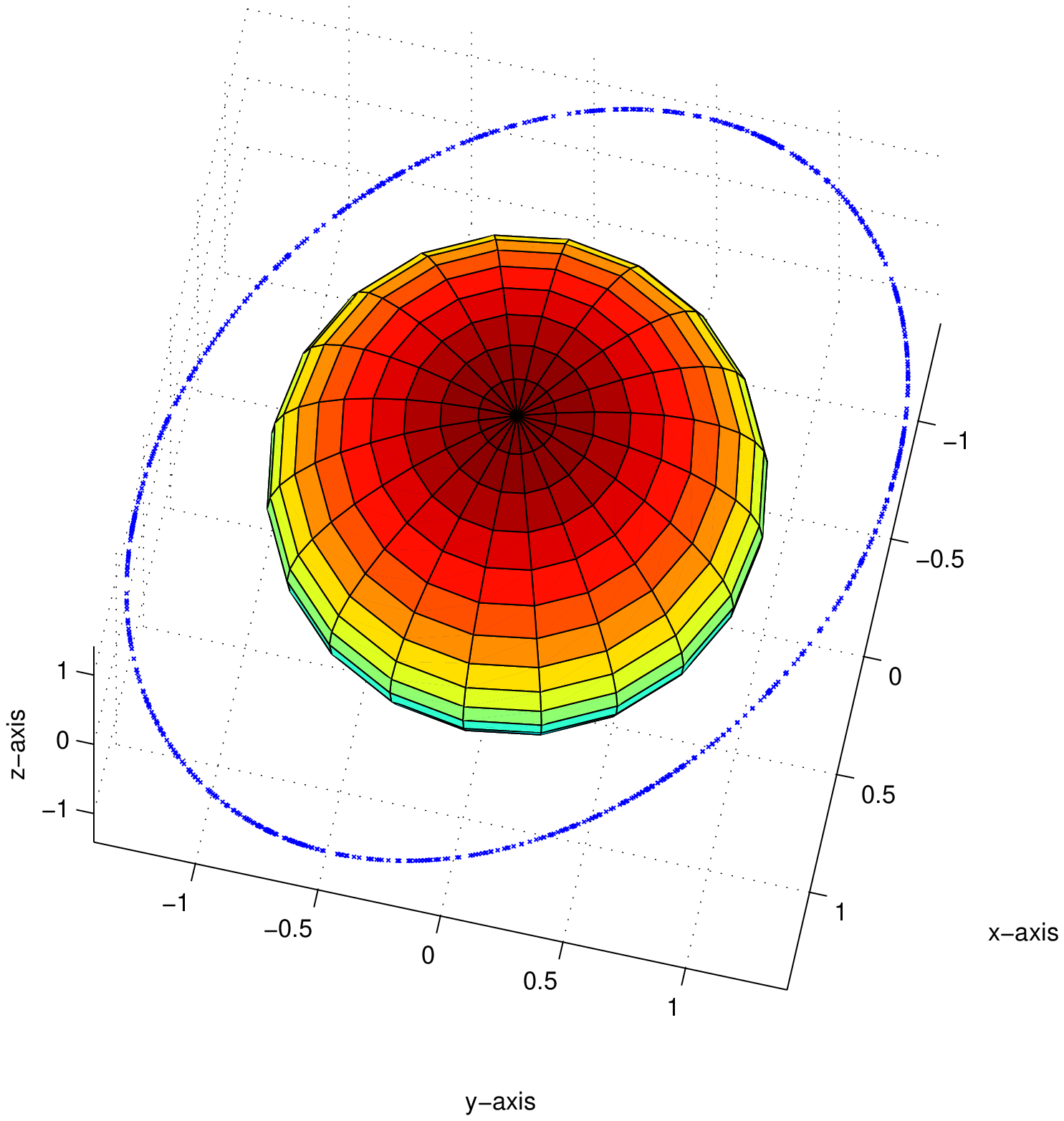}
\caption{Set of fixed points + unit sphere from different perspectives}
\label{fig_sim}
\end{figure}

\newpage

\begin{figure}[!t]
\centering
\includegraphics[width=0.495\textwidth]{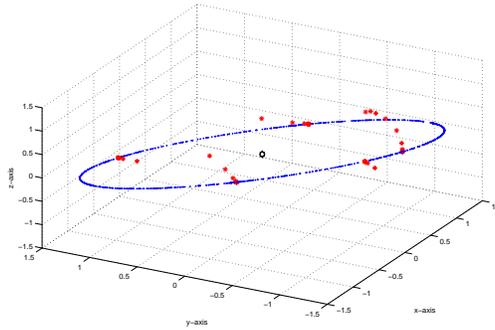}
\caption{Illustration of convergence to fixed points from multiple starts}
\label{fig_sim}
\end{figure}


\begin{figure}[!t]
\centering
\includegraphics[width=0.495\textwidth]{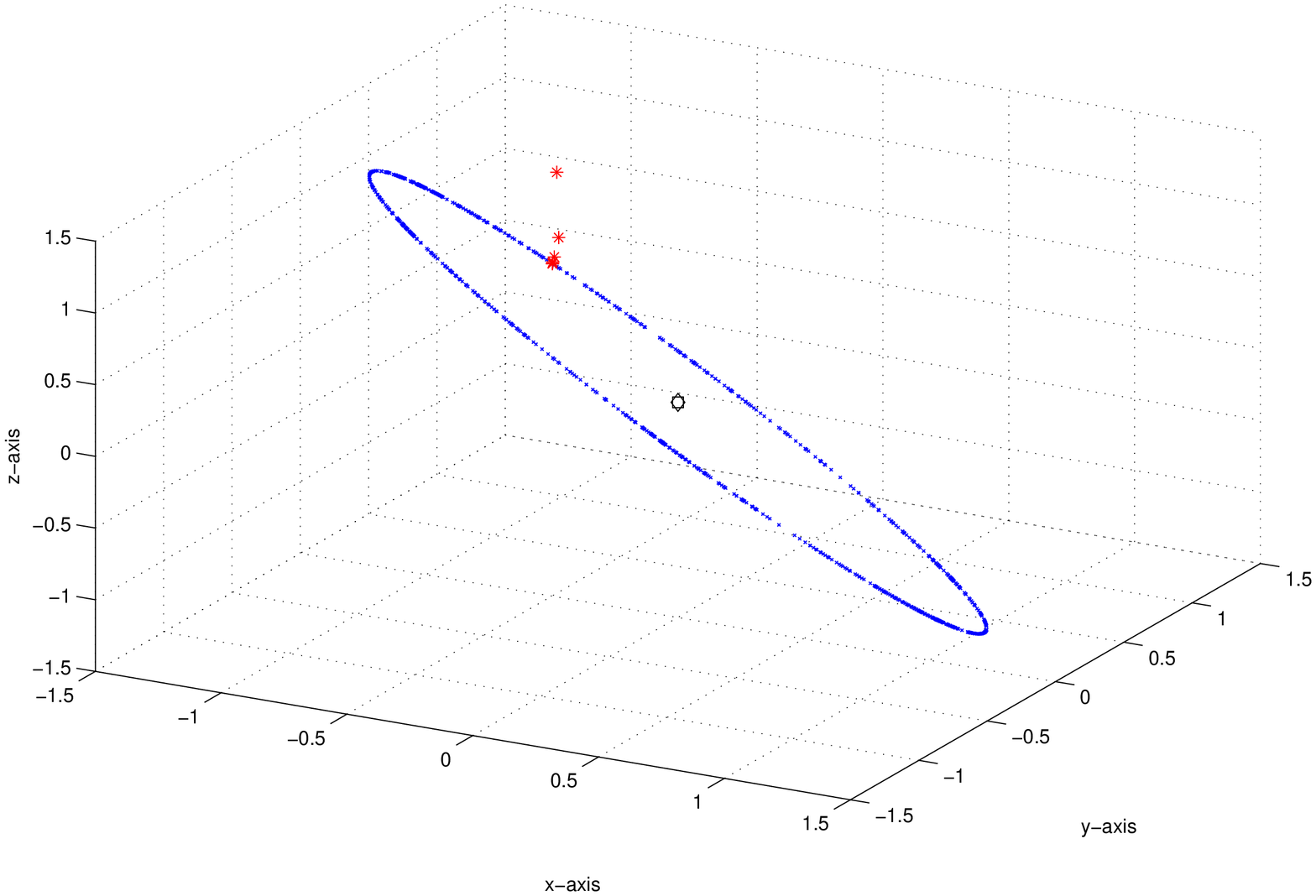}
\includegraphics[width=0.495\textwidth]{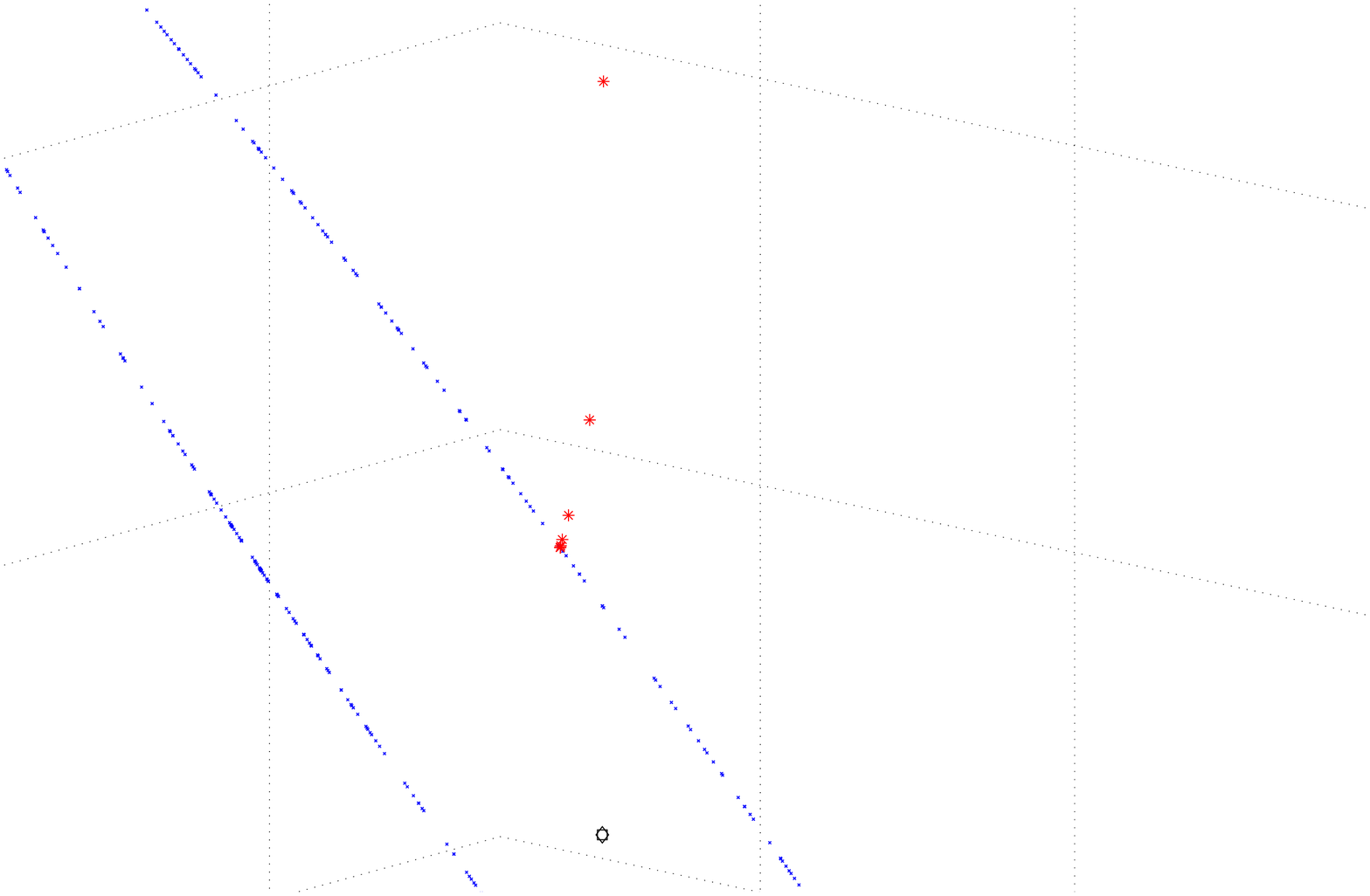}
\caption{(Left) Example of convergence to set of fixed points, (Right) Close up of Example }
\label{fig_sim}
\end{figure}


%


\begin{figure}[!t]
\centering
\includegraphics[width=0.495\textwidth]{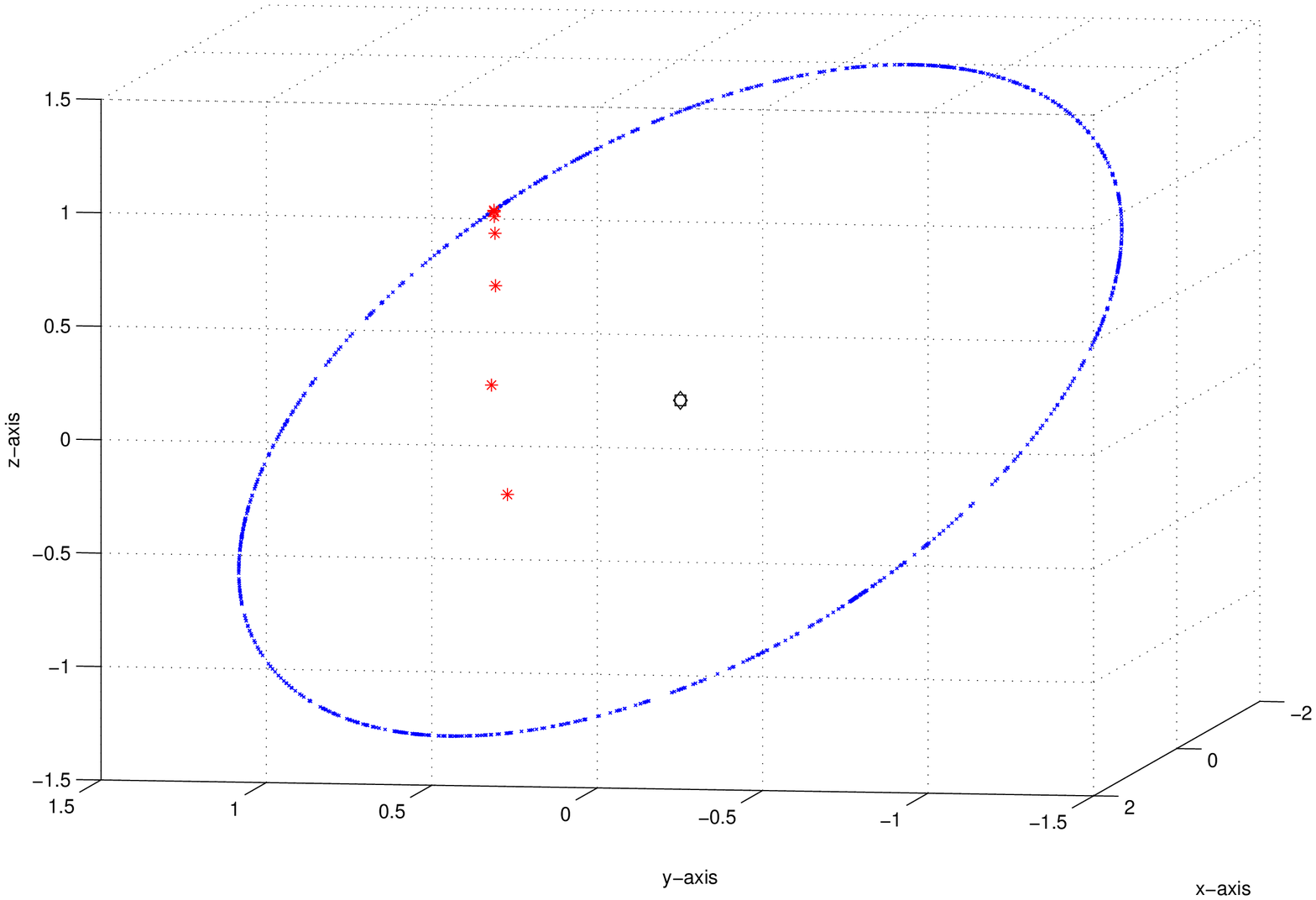}
\includegraphics[width=0.495\textwidth]{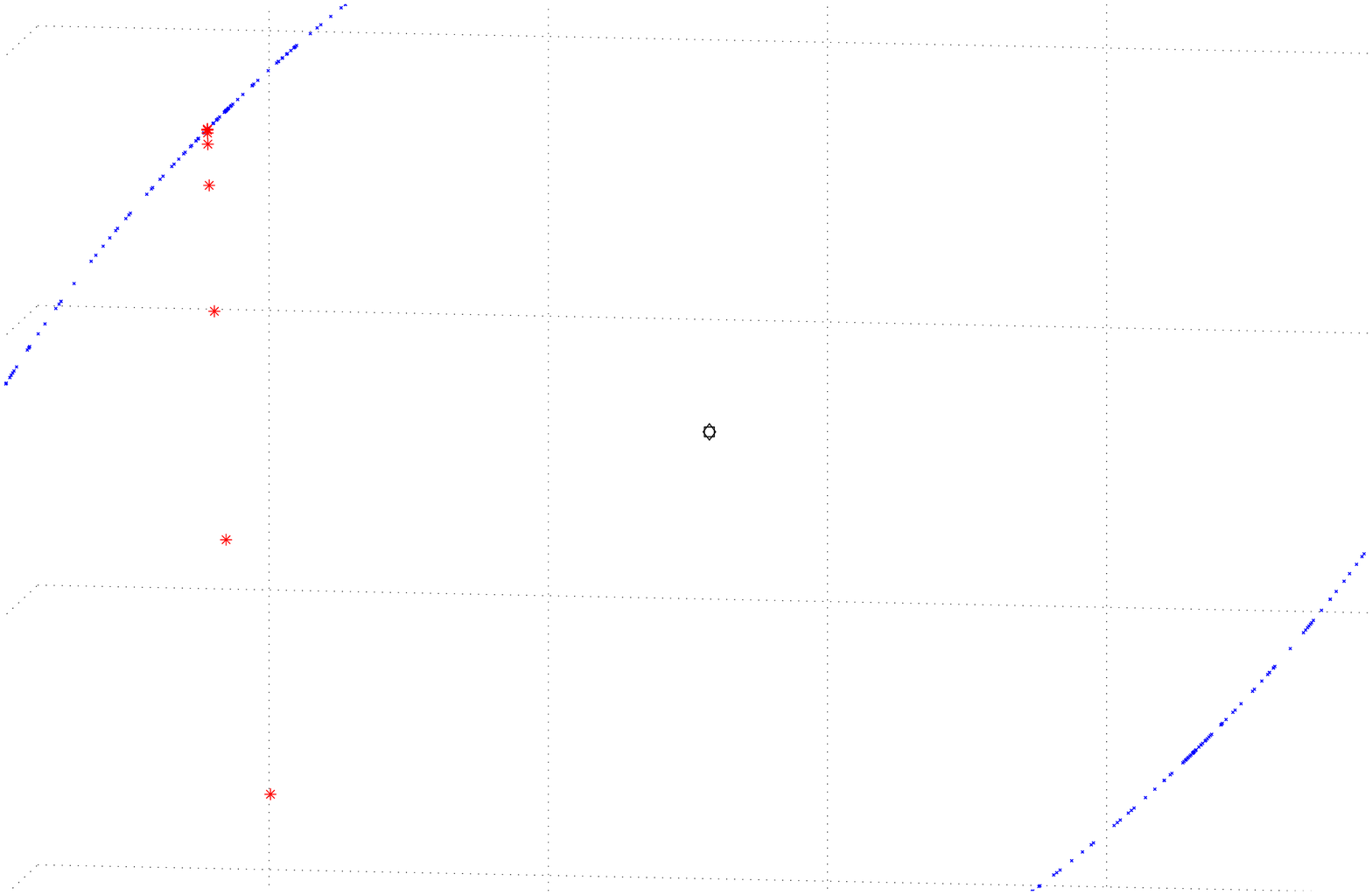}
\caption{(Left) Example of convergence to set of fixed points, (Right) Close up of Example }
\label{fig_sim}
\end{figure}


\begin{figure}[!t]
\centering
\includegraphics[width=0.495\textwidth]{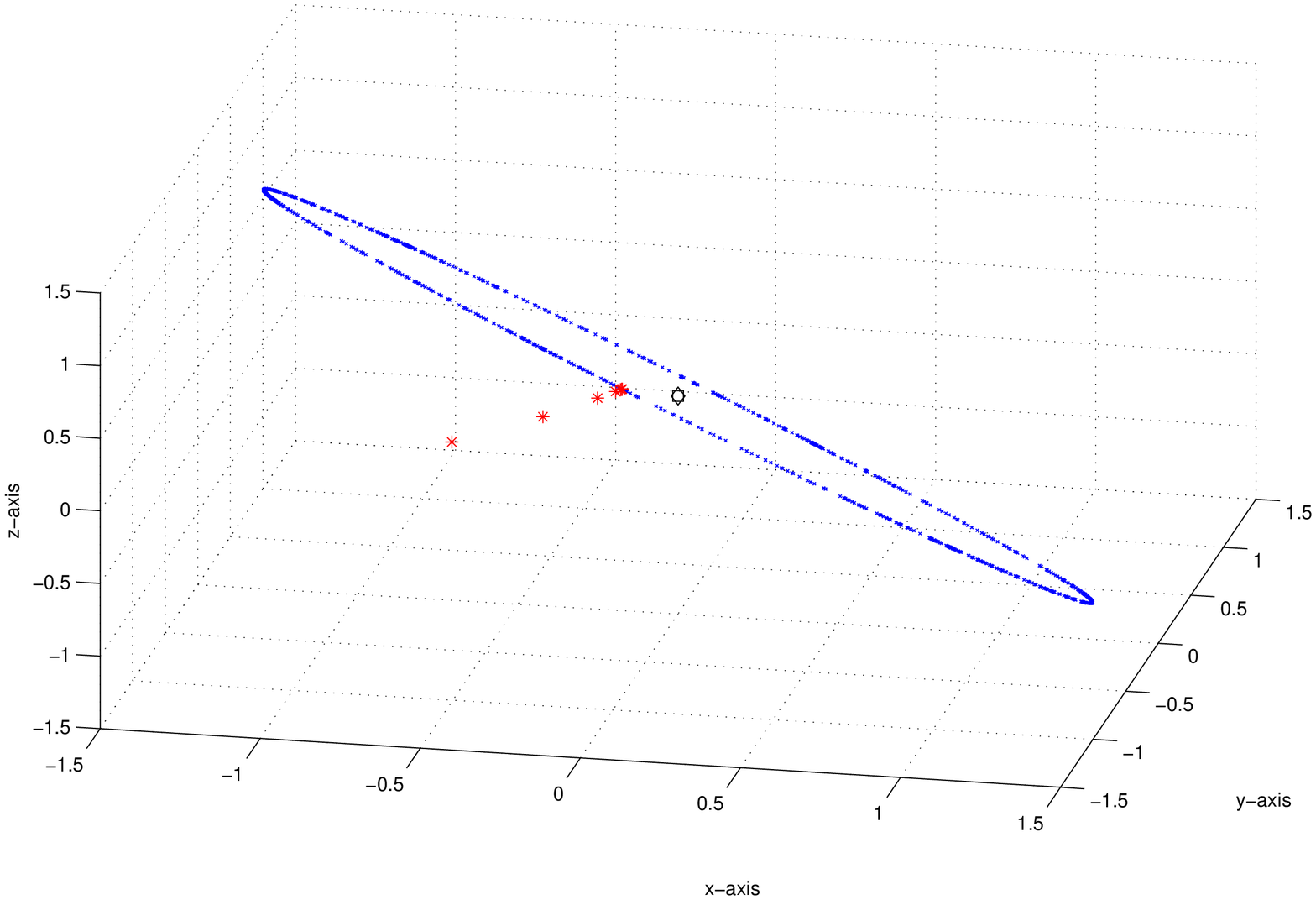}
\includegraphics[width=0.495\textwidth]{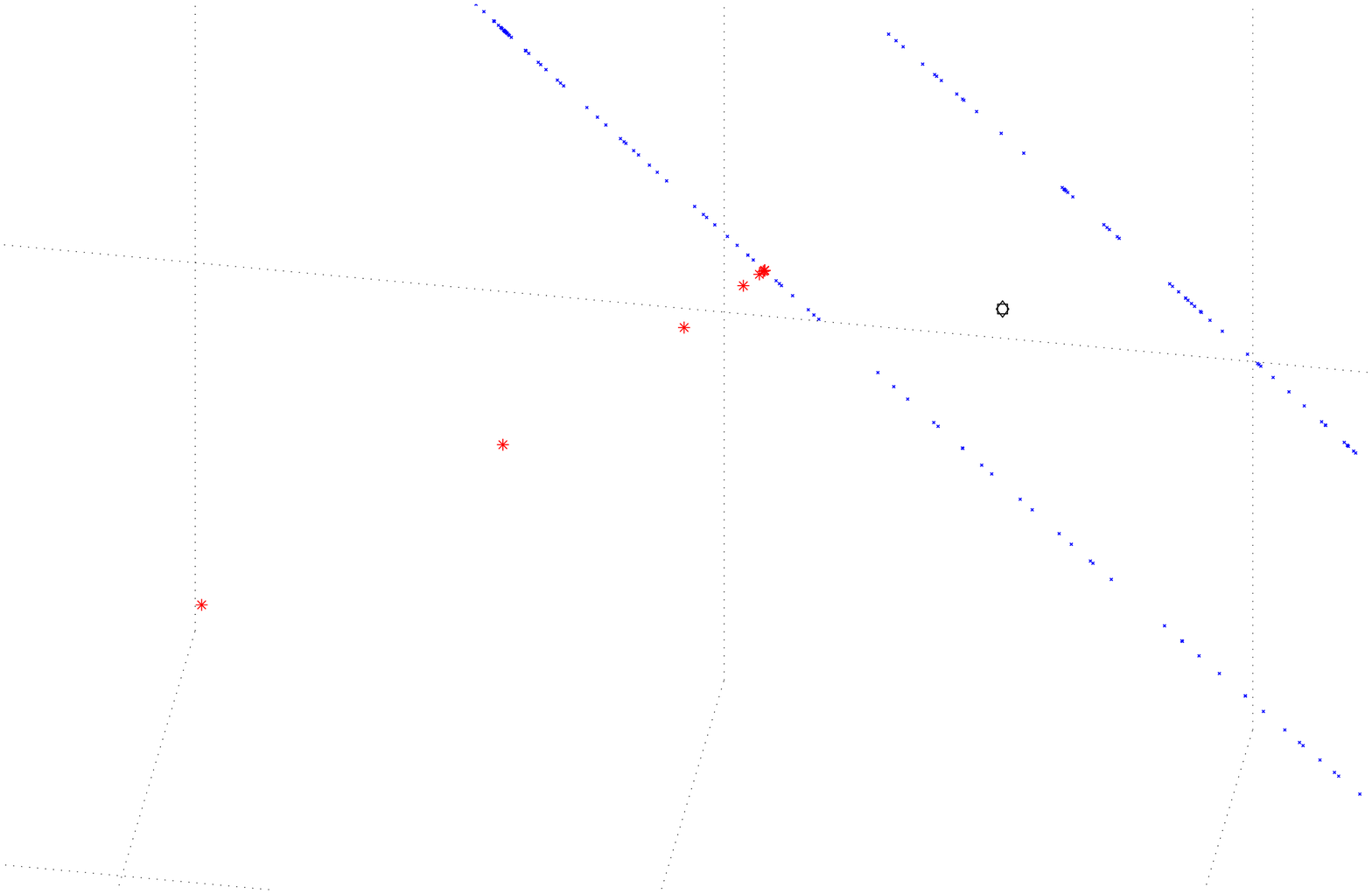}
\caption{(Left) Example of convergence to set of fixed points, (Right) Close up of Example }
\label{fig_sim}
\end{figure}

%
%

\section{Topics for future study}
There are at least several directions for future research that continues to build from \cite{olshen_rajaratnam1}, 
\cite{olshen_rajaratnam2} and from material presented here.  For one, Adam Olshen notes that often in applications, and for many reasons, $\bf{X}$ may not be exactly rectangular. Thus, some initial coordinates may be missing.  By arguments not reported here we have shown that so long as ``missingness" is independent of values that would be reported for a full $\bf{X}$, and no matter the realization there are almost surely genuine data for at least three rows for each column and at least three columns for each row, and so long as standard deviations are defined with ``correct:" divisors, then convergence on all but a set of initial values of measure 0 is plausible. The process of rendering rows and columns with means 0 and standard deviations 1 is potentially a powerful ``hammer" in search of a biological ``nail."  Just as studying data by (scale-free) correlations rather than in original given scales can lend insight, so, too, can the process of successive normalization described here.

\section*{Acknowledgment}

The authors thank colleague Bradley Efron for introducing them to the original problem. They are grateful to David Siegmund, Johann Won and Adam Olshen for useful discussions. David Sigemund told us of a mistake in the argument for Theorem 4.2 that was given in \cite{olshen_rajaratnam1}. It was discovered by him and Michael Hogan and is corrected here. The authors also acknowledge Sang Oh for LaTeX assistance. Richard Olshen was supported in part by grants 4R37EB002784-35 (an NIH MERIT award), 1U19AI090019-01 and UL1 RR025744. Bala Rajaratnam was supported in part by grants NSF-DMS 0906392, NSF-DMS-CMG 1025465, NSF-AGS-1003823, NSA H98230-11-1-0194, DARPA-YFA N66001-11-1-4131, and SUWIEVP10-SUFSC10-SMSCVISG0906.

\bibliographystyle{mdpi}
\makeatletter
\renewcommand\@biblabel[1]{#1. }
\makeatother


\end{document}